\documentclass[12pt]{article}
\usepackage{amsmath,amsfonts,amsthm,amscd}
\newtheorem{theorem}{Theorem}
 
\newtheorem{proposition}{Proposition}
\newtheorem{lemma}{Lemma}

\newtheorem{corollary}{Corollary}
\newtheorem{remark}{Remark}
\numberwithin{equation}{section}
\numberwithin{theorem}{section}
\numberwithin{proposition}{section}
\numberwithin{lemma}{section}
\numberwithin{claim}{section}
\numberwithin{corollary}{section}
\newcommand{\cpn}{\ensuremath{\mathbf{C}P^{N}}}

\newcommand{\dl}{\ensuremath{\partial}}
\newcommand{\dlb}{\ensuremath{\overline{\partial}}}
\newcommand{\gr}{\ensuremath{Gr(N-n,\mathbf{C}^{N+1})}}
\newcommand{\fl}{\ensuremath{\mathbf{Fl}(N-n-1,N-n+1;\mathbf{C}^{N+1})}}
\newcommand{\slnc}{\ensuremath{SL(N+1,\mathbf{C})}}

\newcommand{\bull}{\ensuremath{{}\bullet{}}}
\newcommand{\dmz}{\ensuremath{(N-n)(n+1)-1}}
\newcommand{\zto}{\ensuremath{(N-n)(n+1)-2}}
\newcommand{\gc}{\ensuremath{G^{\mathbb{C}}}}
\newcommand{\Ce}{\ensuremath{\mathbf{C}}}
\newcommand{\ra}{\ensuremath{\longrightarrow}}
\input{diagrams.tex}

\newarrow{to}---->
\newarrow{mto}|--->
\newarrowtail{<=}\Rightarrow\Leftarrow\Uparrow\Downarrow
\newarrow{bboth}{<=}==={=>}
\newarrow{inject}{hooka}---{vee}
\newarrow{surject}----{>>}
\begin{document}
\title{Analysis of Geometric Stability}

\author{{Sean T. Paul}\thanks{Research supported in
part by an NSF Postdoctoral Fellowship} \and
{Gang Tian}\thanks {Research supported in part by an NSF grant and a Simons
fund}}

\date{April 11, 2004}
\maketitle
\vspace{-5mm}
\begin{abstract}
We identify the difference between the CM polarisation
and the Chow polarisation on the ``Hilbert scheme''. 
As a consequence, we give a numerical criterion for the CM stability as in
Mumfords' G.I.T..
Also, we write down an explicit formula for the generalised futaki invariant interms of weights and multiplicities of the associated degeneration.

\end{abstract}
\newpage
\section{Introduction}
The concept of CM stability was introduced by the second
author in early 90's to study the K-energy on the space of
K\"ahler metrics on a compact K\"ahler manifold $X$ with fixed
K\"ahler class. The properness of the K-energy implies the CM
stability of the underlying polarized manifold, this follows from the fact thatthe K-Energy is the logarithm of a \emph{singular} metric on the CM polarisation (cf. [Tian94],
[Tian97]). The purpose of this paper is to relate the CM
stability to Chow points and consequently, the Chow stability
introduced by Mumford. As a corollary of our analysis here, we
also refine the definition of the K-stabiilty first studied in
[Tian97] and establish its numerical criterion. In a related work
Zhiqin Lu [Lu] has analyzed the limit of the derivative of the 
K-Energy map on smooth hypersurfaces in $\cpn$.

We show that the CM stability is equivalent to the stability of what we have called the Double Chow point of $X$. This relates CM stability to conventional
notions from G.I.T. and Chow coordinates. The chow form construction (which is due to Cayley) dates 
back about 150 years, the fact that it is almost always a \emph{singular} divisor
is well known (see [Mum]). Despite this, our construction of the double chow point in this paper seems to be new .
The study of the CM stability and the K-stability was inspired by
the Calabi problem on existence of K\"ahler-Einstein metrics. In
1976, S.T. Yau solved the Calabi conjecture which gave the
existence of K\"ahler-Einstein metrics in the case $c_1(X)=0$.
Around the same time, Aubin and Yau independently proved the
existence of K\"ahler-Einstein metrics in the cases $c_{1}(X)< 0$.
Since then, much effort has been directed towards the remaining
case of positive Ricci curvature. This has proved to be a
difficult and rich area of investigation. One of the main results in
[Tian97] is that the existence of K\"ahler-Einstein metrics is
equivalent to the properness of the K-energy. Hence, the existence
of K\"ahler-Einstein metrics implies the CM stability and the K
stability. In early 80's, S.T. Yau conjectured that the existence of K\"{a}hler
Einstein metrics is related to the asymptotic Chow stability. The necessary half of
Yau's conjecture was affirmed by S. Donaldson in 2001.
He showed that asymptotic Chow stability follows
from the existence of metrics of constant scalar curvature. 
In view of this and results in [Tian97], one is tempted to give a
direct proof that the CM stability implies the Chow stability.
In connection with this see also [Zhang96], [Luo98], and [Paul00]. More recently, X. Wang [Wang02] used these ideas in the context of vector bundles.  
A few years ago, Mueller and Wendland [WenMu99] studied the K-energy in terms of determinant line bundles using the (deep) work of [BGSI,II,III88]. Thier polarisation should be compared with the CM polarisation. Similarly, Phong and Sturm [PhSt02] studied the K-energy in terms of Chow coordinates of underlying manifolds and were led to a certain notion of stability.
It is plausible that this stability is nothing but the CM
stability.

The organization of this paper is as follows. In the next section we give a brief account of the construction of the Chow point and introduce the ``Double'' Chow point. Included here are a few words about Mumfords' G.I.T. and the CM stability. For a fuller account 
we refer the reader to the literature. In section 3 we state our main theorem and some of its immediate corollaries. We also review notions from previous papers that we need.
We defer the proofs to section 4. In section 5 we discuss our refinement of the generalised Futaki Invariant. 

The proof of Theorem 3.1 is an application of GRR and base change, and is completely algebro-geometric.  
Theorem (\ref{asymp}) is essentially Theorem 3.1 carried out at the potential level. Broadly speaking, there are just two steps here. The first step (Theorem (\ref{t94})) consists in expressing stability as an integral over the hypersurface $Z_{X}$, this was carried out in [Tian94]. The second step (Theorem (\ref{p2000})) consists of transforming this integral to an integral over $X$ itself, this step had been taken about three years ago by the first author. All we need to do here is put these two ingredients together.
So our results here are a natural continuation of the work in the first authors dissertation.
The first author would also like to thank his colleagues at the Massachusetts Institute of Technology for his visit during the fall of 2001, where the work in this paper began.
The remaining part of this work was done while the second author was visiting the Department of Mathematics, Columbia University in the spring of 2002. He thanks his colleagues there for providing him with an excellent research
environment.

\section{Chow Forms}
 The idea behind the Chow form is to describe
an arbitrary variety by a \emph{single} equation, this was done by Cayley for curves, then later
generalised by Chow. Since this is our main concern in what follows we give some details. 
For further information we refer the reader to [GKZ].
Let $X \subset \cpn$ be an $n$ dimensional subvariety of $\cpn$ with degree $d$, then the Chow form,
or associated hypersurface to $X$ is defined by
\begin{align}\label{chow}
Z_{X}:= \{L \in \mathbb{G}(N-n-1,\mathbb{C}P^{N}): L\cap X \neq \emptyset\}
\end{align}

It is easy to see that $Z_{X}$ is a hypersurface (of degree $d$) in $\mathbb{G}(N-n-1,\mathbb{C}P^{N})$.
Since the homogeneous coordinate ring of the grassmannian is a UFD, any codimension one subvariety with 
degree $d$ is given by the vanishing of a section
\begin{align}
S \in H^{0}(\mathbb{G}(N-n-1,\mathbb{C}P^{N}),\underline{O}(d)) 
\end{align}
Note that $H^{0}(\mathbb{G}(N-n-1,\mathbb{C}P^{N}),\underline{O}(d))$ is an \emph{irreducible}
$\slnc$ module.\\
In our case there is an 
\begin{align}
R_{X} \in \mathbf{P}H^{0}(\mathbb{G}(N-n-1,\mathbb{C}P^{N}),\underline{O}(d))
\end{align}
such that
\begin{align}
\{\ R_{X}=0\ \}= Z_{X}
\end{align}
By abuse of terminology, we will often call $R_{X}$ the chow form (or chow point) of $X$.
Unfortunately  $Z_{X}$ is almost \emph{never} smooth, the singular set $\mathcal{D}$ is given by
\begin{align}
\mathcal{D}:= \{L \in Z_{X}: \sharp(L\cap X)\geq 2 \}
\end{align}

What's relevant here is  the top strata of this set (also denoted by $\mathcal{D}$)  
\begin{align}
\mathcal{D}:= \{L \in Z_{X}: \sharp(L\cap X)= 2 \}
\end{align}

Since $\mathcal{D}$ has codimension two in
$\gr$ we observe that it has a\\ 
\emph{bidegree} $(d_{1},d_{2})$ which we now describe.
Let $k:= N-n-1$, and fix a flag\\ $E_{k-2}\subset E_{k+1}\in \mathbf{Fl}(k-2, k+1,\mathbb{C}P^{N})$. This flag determines a plane (which we\\ denote by $\sigma_{I}$) in $\gr$
\begin{align}
\sigma_{I}:= \{L\in \gr:E_{k-2}\subset L \subset E_{k+1}\}\approx \mathbb{C}P^2
\end{align}
Similarly, the flag $E_{k-1}\subset E_{k+2}$
determines a plane 
in $\gr$, which we just denote by $\sigma_{II}$, and call them planes of type II.
Then we define
\begin{align} 
d_{1}:= \#(\sigma_{I}\cap \mathcal{D})\ \mbox{and} \ d_{2}:= \#(\sigma_{II}\cap \mathcal{D}).
\end{align}

Of course this makes sense for any codimension two analytic subvariety of\\ 
$\gr$. Below, we calculate the $d_{i}$ in terms of the degree of $X$ in $\cpn$.
The extension Theorem (\ref{ext} below) rests on associating to the singular 
divisor $\mathcal{D}$ in $Z_{X}$ a hypersurface inside of the two step flag variety
\begin{align}
\mathbf{Fl}_{2}:= \fl
\end{align}
in much the same way that $X$ was
associated to $Z_{X}$ . This is the main innovation in the paper. 
We remark that the purpose of the chow form construction is to \emph{minimize}
the codimension, therefore it seems imprudent to first map $\mathcal{V}$ into a projective space via the Pl\"{u}cker embedding (thereby making its codimension large), and then take the chow form of the image in the conventional manner. Many of the features of Cayleys' construction are easily established from our point of view\footnote{Irreducibility of the recieving vector space (as an $\slnc$-module), preservation of degree, and the Bezout formula (see below)}.
Since $\mathcal{D}$ has codimension two in
$\gr$ we need to locate the \emph{lines} in
$\gr$. They are given by \emph{pencils}
determined by an arbitrary point in $\mathbf{Fl}_{2}$, we shall denote such
lines by $P_{E,F}$.
\begin{align}
P_{E,F}:= \{L\in \gr :E\subset L \subset F\}\approx \mathbb{C}P^1
\end{align}
We remark that there is only one kind of line in $\gr$.
Now define $Z_{\mathcal{D}}$ to be all flags $E \subset F$ in 
$\mathbf{Fl}_{2}$ whose associated pencil meets $\mathcal{D}$
\begin{align}
Z_{\mathcal{D}}:=\{E\subset F\in \fl : P_{E,F}\cap\mathcal{D}\neq \emptyset\}
\end{align}
Recall that
\begin{align}
Pic(\fl)=\mathbf{Z}\oplus\mathbf{Z}
 \end{align}
 
The two generators being the determinants of the respective canonical quotient  bundles
\begin{align}
\mathbb{L}_{1}|_{E\subset F}:= Det(\mathbf{C}^{N+1}/E) \qquad \mathbb{L}_{2}|_{
{E\subset F}}:= Det(\mathbf{C}^{N+1}/ F)
\end{align}
It is not difficult to see that  $Z_{\mathcal{D}}$ has codimension 1 in $\mathbf{Fl}_{2}$, hence there  are nonnegative integers $({d_{1}}^{'},{d_{2}}^{'})$ and a unique (up to scale) section 
\begin{align}
f_{\mathcal{D}}  \in H^{0}(\mathbf{Fl}_{2};\ \underline{O}(\mathbb{L}_{1}^{{d_{1}}^{'}}\otimes\mathbb{L}_{2}^{{d_{2}}^{'}}))
\end{align}
such that
\begin{align}
\{f_{\mathcal{D}}=0\}= Z_{\mathcal{D}}
\end{align}
Now, just as  in the classical construction of the chow form we have that 
\begin{align*}
H^{0}(\mathbf{Fl}_{2};\ \underline{O}(\mathbb{L}_{1}^{{d_{1}}^{'}}\otimes\mathbb{L}_{2}^{{d_{2}}^{'}}))
\end{align*}
is also an irreducible module for $\slnc$. 
We will refer to this section $f_{\mathcal{D}}$ simply as the chow point of $\mathcal{D}$.
The association of $f_{\mathcal{D}}$ to $\mathcal{D}$ shares many of the important features of the classical chow point. At least when our subvariety $\mathcal{V}$ has codimension two we have the important fact 
\begin{align*}
\bull \mathcal{V}\ \mbox{is completely determined by}\ f_{\mathcal{V}} 
\end{align*}
In other words\\
\\
\emph{A point} $L\in \gr$ \emph{belongs to} $\mathcal{V}$ 
\emph{if and only if}\emph{ for every} $P_{E,F}$ \emph{containing} $L$ \emph{we have that} $f_{\mathcal{V}}(E,F)=0$.
\\
\\
This amounts to nothing more than the fact that given any point $L$ not on $\mathcal{V}$ there exists a pencil $P_{E,F}$ passing through $L$ which does not meet $\mathcal{V}$.

Recall that for a subvariety $X$ of $\cpn$ the degree of $X$ is the same as the degree of $R_{X}$.
Once again, this holds for our construction.
\begin{proposition}
 $({d_{1}}^{'},{d_{2}}^{'})= (d_{1},d_{2})$
\end{proposition}
{\bf{Proof}}\\
Fix a flag $E_{k-2}\subset E_{k}\subset E_{k+1}$. Then the ``$E$-lines'' in $\mathbf{Fl}_{2}$ are given by
\begin{align}
\ell_{E}:= P_{E_{k-2},E_{k}}\times \{E_{k+1}\}\subset \mathbf{Fl}_{2} 
\end{align}
One defines the ``$F$-lines'' similarly.
We want to count the number of points of intersection of  $Z_{\mathcal{D}}$ with a generic $\ell_{E}$
\begin{align}
\#(\ell_{E}\cap Z_{\mathcal{D}}) 
\end{align}
We know that $\sigma_{I}\cap {\mathcal{D}}=\{L_{1},\dots,L_{d_{1}}\}$, $dim(L_{i})=k$.
Since $E_{k-2}\subset E_{k}\subset E_{k+1}$ is generic, for $i\neq j$ $L_{i}\cap E_{k} \neq L_{j}\cap E_{k}$ and $E_{k}\notin \{L_{1},\dots, L_{d_{1}}\}$. Now it is easy to see that
\begin{align}    
\ell_{E}\cap Z_{\mathcal{D}}=\{L_{i}\cap E_{k},1\leq i \leq d_{1}\}.
\end{align}
Therefore,
\begin{align}
\#(\ell_{E}\cap Z_{\mathcal{D}})= d_{1} .
\end{align}
Similarly one finds 
\begin{align}
\#(\ell_{F}\cap Z_{\mathcal{D}})= d_{2}.
\end{align}
Next we observe that 
\begin{align*}
\underline{O}_{\mathbf{Fl}_{2}}({d_{1}}^{'},{d_{2}}^{'})|\ell_{E} \approx \underline{O}_{\mathbb{P}^{1}}({d_{1}}^{'})
\end{align*}
and $f_{\mathcal{D}}|\ell_{E}$ corresponds to a nonvanishing section of 
$\underline{O}_{\mathbb{P}^{1}}({d_{1}}^{'})$ which has ${d_{1}}^{'}$
zeros, therefore ${d_{1}}^{'}= d_{1}$, the same argument shows that ${d_{2}}^{'}=d_{2}$
${}\Box{}$.
\\
\\
The ``Bezout formula'' holds for $f_{\mathcal{D}}$ as well.  To formulate it let $\mathcal{V}$ and $\mathcal{W}$ be two hypersurfaces in $\mathbf{G}(k,n)$ (with defining functions $F_{1}$ and $F_{2}$ and degrees $d_{1}$ and $d_{2}$ respectively), intersecting transversely. Let $P_{E,F}$ be a generic pencil. Then our Bezout formula is

\begin{align*}
f_{\mathcal{V}\cap \mathcal{W}}(E\subset F)= \prod_{i=1}^{d_{1}}F_{2}(L_{i})
\qquad
P_{E,F}\cap \mathcal{V}= \{L_{1},\dots,L_{d_{1}}\}
\end{align*}

Next we will calculate the bidegree $(d_{1},d_{2})$ of our $f_{\mathcal{D}}$.
It is really enough to consider the case of a curve $E$ in $\mathbb{P}^{3}$. 
Let $d= deg(E)$.
The chow form of $E$, $Z(E)$ lives in $\mathbb{G}(1,\mathbb{P}^{3})$ which has dimension 4, so $\mathcal{D}(E)$
is an algebraic surface. Therefore we look at planes in $\mathbb{G}(1,\mathbb{P}^{3})$ which meet $\mathcal{D}(E)$.
Recall that there are two families of such planes, both parametrised by $\mathbb{P}^{3}$, the first type $\sigma_{1}$ is given by

\begin{align}
\sigma_{1}:= \{L \in \mathbb{G}(1,\mathbb{P}^{3}):L \subset V_{1}, dim(V_{1})=3; V_{1}\subset \mathbf{C}^{4}\}
\end{align}

the second type is
\begin{align}
\sigma_{2}:= \{L \in \mathbb{G}(1,\mathbb{P}^{3}) :\mathcal{U} \subset L , dim(\mathcal{U})=1\}
\end{align}

We just need to calculate
\begin{align}
\#(\sigma_{1}\cap \mathcal{D}(E))= \#\{ l\subset \mathbb{P}(V_{1})\subset \mathbb{P}^{3}: \#(l\cap E)=2 \}
\end{align}

Since we have chosen $V_{1}$ generically we have that $\#(\mathbb{P}(V_{1})\cap E)= d$\\ and the points of intersection are distinct  $ p_{i}\neq p_{j}$ for $i\neq j$. So there are 
\begin{align}
d_{1}=\frac{d(d-1)}{2} 
\end{align}
such lines.

As for $d_{2}$, we look at

\begin{align}
\sigma_{2}\cap\mathcal{D}(E) = \{pt.\in L\in \mathbb{G}(1,\mathbb{P}^{3}): \#(L\cap E)=2\}
\end{align}

Again we can assume that $pt.\notin E$. So we project from this point
\begin{align}
\pi_{pt.}:\mathbb{P}^{3}\rightarrow \mathbb{P}^{2}
\end{align}

and we have
\begin{align}
\#\{L\subset \mathbb{P}^{3}; pt.\in L \thickspace \#(L\cap E)=2\} = \#(nodes \thinspace of\thinspace \pi_{pt.}(E))
\end{align}

Now, this latter number can be computed by simply quoting the degree-genus formula.
So we get that
\begin{align}
d_{2}=\frac{(d-1)(d-2)}{2}-g(E)
\end{align}
where $g(E)$ is the geometric genus of $E$.
Exactly the same reasoning works in higher dimensions, with the degree genus formula replaced by the double point formula.$\Box$
\\
\\
Next we recall some basic notions from Mumford's G.I.T.
For notational\\ ease we just set
\begin{align*}
G^{\mathbb{C}}:= \slnc. 
\end{align*}
Let $v$ be any nonzero vector in a finite dimensional rational $\gc$ module $E$.
We say that $v$ is a \emph{stable} (resp. \emph{semistable}) point for the action of $\gc$ if
the orbit of $v$ is closed and the stabilizer of $v$ in $\gc$ is finite, (resp. if the closure of the orbit does not contain $0$). The concept of stability was introduced by David
Mumford in order to construct moduli of curves, and vector bundles. There is now a vast literature on the subject (see [Mum '82] and references therein). In general, it is still a challenging problem to verify the stability in many concrete cases. The \emph{Numerical criterion}
is perhaps the most important (if not the only) device in the subject. In order to state it, let $\lambda$ be an algebraic 1psg of $\gc$.
Recall that $\mu(\lambda, v)$ is the unique integer (called the \emph{slope}) such that
\begin{align}
\lim_{\alpha \rightarrow 0}\alpha^{\mu(\lambda, v)}\lambda (\alpha)v
\end{align} 
exists and is not equal to $0$.
Then we have the following fundamental result\\
\begin{center}(\emph{Hilbert Mumford Numerical criterion})
\end{center}
\begin{align*}
\bull &v \ \mbox{is stable iff}\ \mu(\lambda, v)>0\ \mbox{for every 1psg}\ \lambda \\ \\
\bull &v \ \mbox{is semistable iff}\ \mu(\lambda, v)\geq 0\ \mbox{for every 1psg}\ \lambda 
\end{align*}
Since we do use the slope in the context of $\gc$-linearisations, we say a few words about it. A line bundle $\mathcal{L}$ over a $\gc$ manifold $\mathcal{H}$ is said to be \emph{linearised} wrt $\gc$ if there is a lift of the given action of $\gc$ on $\mathcal{H}$ 
to  $\mathcal{L}$ compatible with the projection
\begin{align*}
\sigma :\mathcal{L}_{x}\rightarrow \mathcal{L}_{\sigma x}\quad \sigma \in \gc
\end{align*}
Assume that $\mathcal{H}$ is compact, then any $x \in \mathcal{H}$ has a limit
$x_{\infty}$ under the action of a 1psg $\lambda$ of $\gc$.
Clearly $x_{\infty}$ is fixed by $\lambda$, this gives a 1 dimensional representation 
$\mathcal{L}_{\infty}$ of $\mathbb{C}^{*}$ via $\lambda$ and hence a character. 
This is $-\mu$.
\begin{align*}
\lambda(\alpha)v= \alpha^{-\mu(\lambda,v)}v \quad v \in \mathcal{L}_{\infty}
\end{align*}
Of course, in the definition of $\mu$ scaling $v$ has no effect.
In our application, $\mathcal{L}$ is the CM-polarisation.
With these preliminary remarks, we can introduce the main concept in the paper.\\
{\bf{Definition}}\\
Let $X$ be an irreducible $n$ dimensional subvariety of $\cpn$, $\mu(X)$ denotes the average of the scalar curvature of $\omega_{FS}$. Let $R_{X}$ be the chow point of $X$. We define the \emph{Double Chow point} of $X$ to be the pair
\begin{align}
([R_{X}],\ [f_{\mathcal{D}}]) \in\mathbf{P}H^{0}(\gr,\ \underline{O}(d))\times
\mathbf{P}H^{0}(\mathbf{Fl}_{2};\ \underline{O}(\mathbb{L}_{1}^{d_{1}}\otimes\mathbb{L}_{2}^{d_{2}}))
\end{align} 
We say that the Double chow point is \emph{stable} iff
\begin{align}
\left(2d +\frac{\mu(X)}{n+1}-(n+2)\right)\mu(\lambda,\ R_{X})- \mu(\lambda,\ f_{\mathcal{D}})>0
\end{align} 
For every algebraic 1psg $\lambda$ of $\gc$.

{\bf{The Hilbert Scheme}}\\
Now we consider $X$ as a generic (smooth) member $X_{z}$ of a family. Let
\[ \begin{CD}
\mathfrak{X} @>>> \cpn \\
@V\pi_{1}VV \\
\mathcal{H}\\
\end{CD}
\]
be a $G^{\mathbb{C}}$ equivariant holomorphic fibration between smooth projective varieties satisfying
\\
\\
 \begin{tabular}{ll}
 $\bull$ $\mathfrak{X} \subset \mathcal{H}\times \cpn$ is a family of subvarieties of 
 generic dimension $n$,
(i.e. some of the\\ fibers are allowed to have dimension higher than $n$) where the action
\\ of \gc  on $\mathfrak{X}$ is
 induced by the standard action on \cpn.
 \end{tabular}
 \\
 \\
 On $\mathfrak{X}$ consider the virtual bundle
 \begin{align}
 \mathcal{E}:= (n+1)(\mathcal{K}^{-1}-\mathcal{K})\otimes(L-L^{-1})^{n}-n\mu(L-L^{-1})^{n+1}
 \end{align}
 Here $\mathcal{K}$ is the relative canonical bundle. Let $\mathcal{L}_{\mathcal{H}}$ be the inverse of the determinant of the direct image
 of $\mathcal{E}$. We call $\mathcal{L}_{\mathcal{H}}$ the CM polarisation. Then \gc   acts naturally on $\mathcal{L}_{\mathcal{H}}^{-1}$. A point $z\in \mathcal{H}_{\infty}$ (the smooth locus of $\mathcal{H}$ )  is said to be
 \emph{CM stable} wrt $\mathcal{L}_{\mathcal{H}}$, if the orbit of any lift of $z$ to $\mathcal{L}_{\mathcal{H}}^{-1}$
 is closed and the stabilizer is finite, $z$ is \emph{semistable} if the closure of the orbit of any lift of $z$ does not meet the zero section.
\section{Statement of Results}
We set
\begin{align*}
&\nu_{1}:= (n+1)(2d + \frac{\mu(X)}{n+1}-(n+2))\\
&\nu_{2}:= \frac{(n+1)}{D}\\
&D:= \mbox{deg}\ {\gr}
\end{align*}

The central result in this paper is the following Theorem
which identifies the CM polarisation.
\begin{theorem}\label{ext}
\end{theorem}
\emph{Assume that the two maps are injective}\footnote{This will be true in most cases, in fact the first map will be injective for any scheme connected and proper over $\mathbb{C}$}
\begin{align}
&\alpha:\mbox{\mbox{Pic}}^{\gc}(\mathcal{H})\rightarrow \mbox{Pic}(\mathcal{H})\\
&c_{1}:\mbox{Pic}(\mathcal{H})\rightarrow \mbox{H}^{2}(\mathcal{H},\mathbb{Z})
\end{align}
\emph{Then we have the following isomorphism of line bundles on} $\mathcal{H}$.
\begin{align}
L_{CM}^{-1}= \iota^{*}_{HC}\underline{O}_{{\mbox{Div}}_{(d,\dots,d)}}(-\nu_{1})\otimes {L}_{\Delta}^{\nu_{2}}
\end{align}

Where we have defined $\iota_{\mbox{HC}}$ to be the \emph{Hilbert Chow}\footnote{For a thorough account of these matters, see Kollar [Kol95]}
morphism
\begin{align*}
\iota_{\mbox{HC}}:\mathcal{H}\rightarrow {\mbox{Div}}_{(d,\dots,d)}\subset B_{d}\end{align*}
${\mbox{Div}}_{(d,\dots,d)}$ denotes the chow variety of dimension n degree d cycles in $\cpn$.
${L}_{\Delta}$ is the invertible sheaf on $\mathcal{H}$ coming from the singularities of the chow divisor. Details appear below.

As a consequence of this, we deduce the following weight identity ($w_{\lambda}$ is used to denote the weight of the action on the indicated line bundle).
\begin{theorem}\label{weight}For every algebraic 1psg $\lambda$ of $\gc$ we have
\begin{align}
w_{\lambda}(L_{CM}^{-1},z)= \nu_{1}w_{\lambda}(Chow(z))+\nu_{2}w_{\lambda}({L}_{\Delta},z) 
\end{align}
\end{theorem}
When $X_{z}^{\lambda(0)}$ is reduced, we can be more explicit about the weight of the $\mathbb{C}^{*}$ action on the line ${L}_{\Delta}$.
\begin{proposition}Assume $X_{z}^{\lambda(0)}$ has no multiple components then
\begin{align}
w_{\lambda}({L}_{\Delta},z) = -Dw_{\lambda}(f_{\mathcal{D}_{z}}) 
\end{align}
\end{proposition}

This yields the following corollary 
\begin{corollary}Assume that the limit cycle $X^{\lambda(0)}$ is reduced, then 
\begin{align}
w_{\lambda}(L_{CM}^{-1},z)= \nu_{1}w_{\lambda}(Chow(z))-(n+1)w_{\lambda}(f_{\mathcal{D}_{z}})
\end{align}
\end{corollary}
On the potential level, we prove the following.
\begin{theorem}\label{asymp}(Asymptotics of the K-energy map)\\
Let $X\subset \cpn$ then\\ 
\begin{align}
d\nu_{\omega}(\varphi_{\sigma})= (2d +\frac{\mu(X)}{n+1}-(n+2))\frac{\log||R_{X}\circ \sigma^{-1}||^{2}}{||R_{X}||^{2}} - \frac{\log||(f_{\mathcal{D}}\circ \sigma^{-1}) ||^{2}}{||f_{\mathcal{D}}||^{2}}+ \frac{1}{D}\Psi_{B}(\sigma)
\end{align}
\end{theorem}
The term $\Psi_{B}(\sigma)$ degenerates to $-\infty$ when the limit cycle has a component of multiplicity greater than one.
\begin{remark}
In view of this, one is tempted to try and prove that
\begin{align*}
\mu(\lambda; R_{X}) \leq  \mu(\lambda;f_{\mathcal{D}}).
\end{align*}
Our bidegree computation lends some plausibility to this.
This would show that K-stability implies Chow (hence Hilbert) stability.
The existence of the Gieseker Scheme follows immediately. However, at the moment, this inequality seems out of reach.
\end{remark}

Very recently, S.K. Donaldson [SKD02] has proposed the following definition of the K-stability, which involves G.I.T. \emph{directly}. To explain this recall the construction of the mth Hilbert point of a subscheme of $\cpn$. So 
let $X$ be a projective variety (or scheme) in $ \cpn$. Then we have, for
$m>>0$, the exact sequence

\[ \begin{CD}
0@>>>H^{0}(\mathcal{I}_{X}(m))@>i>>H^{0}(\cpn,  {\mathcal
O}(m))@>Res_{m}(X)>> H^{0}(X,  {\mathcal O}(m)_{X}) @>>>0
\end{CD}
\]
where the vector space on the left is
\begin{align*}
\begin{split}
H^{0}(\mathcal{I}_{X}(m)))= &\{\mbox{All homogeneous polynomials of degree }\ m\\
&\mbox{in}\ N+1 \ \mbox{variables that vanish on}\ X\}
\end{split}
\end{align*}
If we let $P(m)$ be the Hilbert polynomial of $X$, then we have
\begin{align*}
\begin{split}
&H^{0}(\mathcal{I}_{X}(m))\in Gr(P(m),H^{0}(\cpn,  {\mathcal O}(m))) \\ \\
\end{split}
\end{align*}
where $Gr(P(m),H^{0}(\cpn,  {\mathcal O}(m)))$ denotes the
Grassmannian of \emph{codimension} $P(m)$ subspaces of
$H^{0}(\cpn,  {\mathcal O}(m))$. Using the Plucker embedding, we
may associate to $X\subset \cpn$ the point in the following
projective space
\begin{align*}
 {\bf{det}}(H^{0}(\mathcal{I}_{X}(m)))\in
 \mathbb{P}(\bigwedge^{\binom{N+m}{m}-P(m)}H^{0}(\cpn,  {\mathcal O}(m))).
\end{align*}
The $\mbox{m}^{th}$ Hilbert point of $X$ with respect to given
polarization ${\mathcal O}(1)|_X$ is given by its dual
${\bf{det}}(H^{0}(X,  {\mathcal O}(m)_{X})$. To fix notation, we
will denote this Hilbert point by
\begin{align*}
\mbox{Hilb}_{m}(X):= {\bf{det}}(H^{0}(X,  {\mathcal O}(m)_{X}) \in
\mathbb{P}(\bigwedge^{P(m)}H^{0}(\cpn,  {\mathcal O}(m))^{*}).
\end{align*}
Since $\gc$ acts on this big projective space, we can associate a
weight to each 1psg. $\lambda:\mathbb{C}^{*}\rightarrow \gc$, namely the weight\footnote{recall that the \emph{weight} ($w_{\lambda}$) and the \emph{slope} ($\mu$) differ by a -1} the action on this point. The 1psg $\lambda$ induces an action on
$E=\bigwedge^{P(m)}H^{0}(\cpn,  {\mathcal O}(m))^{*}$. It is easy
to see that any such an action can be diagonalized on $E$. Let
$\{e_{1},\dots,e_{d}\}$ be such a basis, i.e.
$\lambda(\alpha)e_{i}= \alpha^{m_{i}}e_{i} \quad \left(m_{i}\in
\mathbb{Z}\right)$. Next express any $v\in E$ in terms of this
basis $v =\sum_{i=1}^{d}v_{i}e_{i} $. Then the \emph{slope} of
$\mbox{Hilb}_{m}(X)$ relative to $\lambda$ is the number (usually
denoted by $\mu(\lambda,v)$:
\begin{align}
\mbox{Max}\{-m_{i}|v_{i}\neq 0\}
\end{align}
We define the weight $w(\lambda, \mbox{Hilb}_{m}(X))$ of
$\mbox{Hilb}_{m}(X)$ to be the weight $\mu(\lambda,v)$ for any $v$
lifting ${\bf{det}}(H^{0}(X,  {\mathcal O}(m)_{X})$.

Our aim is to study the weight $w({\lambda}, \mbox{Hilb}_{m}(X))$.
By general nonsense, this weight is given by a numerical
polynomial of degree at most $n+1$ where $n=\mbox{dim}(X)$. In
other words
\begin{align*}
w_{\lambda}(\mbox{Hilb}_{m}(X))= a_{n+1}m^{n+1}+
a_{n}m^{n}+O(m^{n-1})
\end{align*}

Also, as is very well known we have, for $m>>0$
\begin{align*}
mP(m)= mh^{0}(X, {\mathcal O}(m))=
b_{n+1}m^{n+1}+b_{n}m^{n}+O(m^{n})
\end{align*}
Where the $b_{i}$ are given by Hirzebruch Riemann-Roch. Following
Donaldson, we let $F_{1}$ be the coefficient of $\frac{1}{m}$ in
the expansion below
\begin{align*}
\begin{split}
\frac{w_{\lambda}(Hilb_{m}(X))}{mP(m)}&= \frac{a_{n+1}m^{n+1}+ a_{n}m^{n}+O(m^{n-1})}{b_{n+1}m^{n+1}+m^{n}+O(m^{n})}\\
&= \frac{a_{n+1}}{b_{n+1}}+
\frac{a_{n}b_{n+1}-a_{n+1}b_{n}}{b_{n+1}^{2}}\frac{1}{m}+
O(\frac{1}{m^{2}})
\end{split}
\end{align*}
A simple computation shows that

\begin{align*}
F_{1}(\lambda, X):= \frac{n!}{2d}(2a_{n}-\mu a_{n+1})
\end{align*}

\emph{In [Don02] S.Donaldson defines  the generalised Futaki
invariant of the degeneration $\lambda$ to be this $F_{1}$.}
\\
\\

In that paper he observes that this coincides with the definition
of Tian \emph{when the central fiber is smooth}. Very recently the authors of this paper have refined this observation, the details will appear in a forthcoming note. Here we introduce what we call the \emph{\bf{Reduced K-Energy map}}.
\begin{theorem} (\emph{Asymptotics of the reduced K-energy map})
 There is a function
$\Psi_{X}:\gc \rightarrow \mathbb{R}$ depending only on the
embedding of $X$ where $-\infty \leq \Psi_{X} \leq C$ such that
\begin{align}
 d\nu_{\omega}(\varphi_{\lambda(t)})-\Psi_{X}({\lambda(t)})= 4 d
 F_{1}(\lambda;X)\log(t)+O(1),
\end{align}
where $\nu_\omega$ denotes the K-energy of Mabuchi.
\end{theorem}
The function  $\Psi_{X}$ is explicit and
degenerates to $-\infty$ if $X^{\lambda(0)}=\lim_{t\mapsto 0}
\lambda(t) (X)$ is non-reduced. The \emph{Reduced} K-energy is
defined to be the quantity on the left hand side of the above
equation. The crucial point is that $\Psi_X$ is \emph{bounded from
above}. This term appears when one compares the CM-stability with
the extremal behavior of the K-energy map. Note that when the limit cycle
$X^{\lambda(0)}$ is smooth we just get Donaldsons' result. The upshot is that these two stability notions are the same- they are just the weight of the CM polarisation introduced by the second author about 10 years ago.
The asymptotics of the (ordinary) K-energy map for an \emph{arbitrary} 
$\lambda$ do not seem to be a genuine G.I.T. notion. This has been pointed out to us by Zhiqin Lu, see [Lu01].

\begin{remark}
 To $X$ we associate not one, but two $\gc$ modules with distinguished points
\\
\\
\begin{tabular}{ll}
$\bull$ $R_{X}\in B_{d}:= \mathbb{P}H^{0}(\mathbb{G}(N-n-1,\mathbb{C}P^{N}),\underline{O}(d))$
\\
\\
\centerline{and}
\\
\\
$\bull$ $f_{\mathcal{D}} \in B_{d_{1},d_{2}}:= \mathbb{P}H^{0}(\fl;\mathbb{L}_{1}^{d_{1}}\otimes\mathbb{L}_{2}^{d_{2}})$
\end{tabular}
\\
\\
Our results give a  relationship between these two points.
This is interesting since the two representations are irreducible and
distinct.
\end{remark}

Lets now go into precisely what we will require. 
Let $\varphi$ be a K\"ahler potential on some $(X,\omega_{FS})\subset \cpn$, and $V:=\int_{X} \omega^{n}$.
Recall the following functionals \\
\begin{align}
\begin{split}
&\bull \qquad J_{\omega}(\varphi):= \frac{1}{V}\int_{X}\sum_{i=0}^{n-1}\frac{n-i}{n+1}
\dl \varphi \wedge \dlb \varphi \wedge \omega^{i}\wedge (\omega + \dl \dlb \varphi )^{n-i-1}\\
 &\bull \qquad F_{\omega}^{0}(\varphi):= J_{\omega}(\varphi)-\frac{1}{V}\int_{X}
\varphi \omega^{n}
\end{split}
\end{align}
Let $\varphi_{t}$ is a smooth path in $ P(X,\omega_{FS})$ (all K\"ahler potentials) joining $0$ with $\varphi$.
Then the \emph{K-energy map}, introduced by Bando and Mabuchi [BM] is given by
\begin{align}
\bull \qquad \nu_{\omega}(\varphi):= -\frac{1}{V}\int_{0}^{1}\int_{X}\dot{\varphi_{t}}(Scal(\varphi_{t})-\mu)\omega_{t}^{n}dt
\end{align}

$Scal(\varphi_{t})$ denotes the scalar curvature of the metric $\omega +\sqrt{-1}\dl\dlb\varphi_{t}$.
It is well known that the K-energy map can be written in the form
\begin{align}
\begin{split}
\bull &\nu_{\omega}(\varphi)=  \frac{1}{V}\int_{M}\log\left(
\frac{\omega_{\varphi}^{n}}{\omega^{n}}\right)\omega_{\varphi}^{n}-\frac{1}{V}\int_{M}h_{\omega}(\omega^{n}- \omega_{\varphi}^{n}) \\
&- \frac{\sqrt{-1}}{V}\sum_{i=0}^{n-1}\int_{M}\dl \varphi \wedge \dlb \varphi \wedge \omega^{i} \wedge \omega_{\varphi}^{n-i-1}
\end{split}
\end{align}

As usual, we define $\varphi_{\sigma}\in P(X,\omega_{FS})$ by the equation
\begin{align}
\sigma^{*}\omega_{FS}= \omega_{FS}+\sqrt{-1}\dl\dlb\varphi_{\sigma}. 
\end{align}

The following Proposition is immediate

\begin{proposition}
There exists an integer $m$ and a constant $C$ such that
\begin{align}
\exp( \nu_{\omega}(\varphi_{\sigma}))\leq C \left(\sum_{i,j}|\sigma_{ij}|^{2}\right)^{m}
\end{align}
\end{proposition}
We will use this below.

Recall the following
\begin{theorem}([Zhang96] [Paul00])\label{p2000}\\
There is a continuous metric $|| \ ||_{u}$ on $\underline{O}_{B}(-1)$ such that
\begin{align}\label{fo}
-\log \left(\frac{{||R_{X}\circ \sigma^{-1}||}^{2}}{{||R_{X}||}^{2}}\right)=
 (n+1)VF_{\omega}^{0}(\varphi_{\sigma})
\end{align}
\end{theorem}

Let $Z_{f}:= \{f=0\}$ be a normal hypersurface of degree $d>1$ in $\mathbb{C}P^{n+1}$ and set
\begin{align*}
B:= \mathbf{P}H^{0}(\mathbb{C}P^{n+1},\underline{O}(d))
\end{align*}
In the next result $\nu_{\omega}$ denotes the K-Energy on the hypersurface $Z_{f}$ and 
$\underline{O}_{B}(-1)$ is the universal bundle over the projective space $B$.

\begin{theorem}([Tian94])\label{t94}
There is a singular hermitian metric $|| \ ||_{Q}$ on $\underline{O}_{B}(-1)$
such that
\begin{align}\label{hyp}
\nu_{\omega}(\varphi_{\sigma})=
\frac{(d-1)(n+2)}{n+1}\log\left(\frac{|| \ ||_{Q}^{2}(f\circ \sigma^{-1})}{|| \ ||_{Q}^{2}(f)}\right)
\end{align}
\end{theorem}
\begin{tabular}{ll}
\emph{Discussion}
\end{tabular}

\begin{tabular}{ll}
$\bull$ The identity (\ref{fo}) says that the stability of $R_{X}$ under the action of $G^{\mathbb{C}}$\\
is completely dictated by the behavior of the energy at $\infty$ in $G^{\mathbb{C}}$.
This holds\\ 
because the norm $|| \ ||$
is continuous (it is actually an  \emph{extension} of the Quillen\\ metric), hence there are positive constants $C_{1}$ and $C_{2}$ such that
\end{tabular}
\begin{center}
$C_{1}||\sigma R_{X}||_{\underline{O}_{B}(-1)}\leq ||\sigma R_{X}|| \leq C_{2}||\sigma R_{X}||_{\underline{O}_{B}(-1)}$ 
\end{center}

\begin{tabular}{ll}
$\bull$ The identity (\ref{hyp}) is concerned with ``how bad''the singularities on a\\ 
hypersurface are allowed to be before the hypersurface destabilizes, this turns out\\ to be quite a difficult 
and delicate issue. Below (Proposition \ref{doub}) we show that\\ 
(\ref{hyp}) holds (with modifications by $f_{\mathcal{D}}$) in the general case.\\
A consequence of (\ref{hyp}) is the following\\
\end{tabular}

\begin{theorem}([Tian94])
{If} $Z_{f}$ {admits a K\"{a}hler Einstein orbifold metric then} $f$ {is semistable}.
\end{theorem}
We have the following, a proof of proposition \ref{fd} will appear in a later section, for the others, we refer the reader to [Paul] and [Tian].
Below $\varphi$ denotes a test form of appropriate type on $G^{\mathbb{C}}$.\\
\\
For a subvariety $\mathcal{V}\subset \gr$ we set
\begin{align*}
G^{\mathbb{C}}\mathcal{V}:= \{(\sigma,L)\in G^{\mathbb{C}}\times \gr: L \in \sigma
\mathcal{V}\}
\end{align*}

\begin{proposition}\label{zx}
Let $X^{n} \subset \cpn$ and let $Z_{X}$ be the associated hypersurface (\ref{chow}), then
\begin{align}
\int_{G^{\mathbb{C}}Z_{X}}
p_{2}^{*}(\omega_{P})^{(N-n)(n+1)}\wedge \pi^{*} \varphi =
D\int_{G^{\mathbb{C}}}\dl\dlb\log{ ||R_{X}\circ \sigma^{-1}||}^{2}
\wedge \varphi
\end{align}
\end{proposition}     

\begin{proposition}\label{fd}
Let $\mathcal{V}$ be a codimension two analytic subvariety of $G(k,\cpn)$
with bidegree $(d_{1}, d_{2})$ and let $f_{\mathcal{V}}$ be its chow from. Then
\begin{align}
\int_{G^{\mathbb{C}}\mathcal{V}}p_{2}^{*}\omega_{Pl}^{(N-k)(k+1)-1}\wedge p_{1}^{*}\varphi
= D \int_{G^{\mathbb{C}}}\dl\dlb \log||(f_{\mathcal{V}}\circ \sigma^{-1}) ||^{2}\wedge 
\varphi
\end{align}
\end{proposition}

\begin{proposition}\label{ken} (complex hessian of the K-energy map)\\
Let $X\subset \cpn$ be a normal subvariety, then
\begin{align}
d\int_{G^{\mathbb{C}}}\nu_{\omega}\dl\dlb \varphi = \int_{G^{\mathbb{C}}X}(R_{G^{\mathbb{C}}|X}+\frac{\mu(X)}{n+1}p_{2}^{*}(\omega_{FS}))\wedge p_{2}^{*}(\omega_{FS}^{n})\wedge p_{1}^{*}\varphi
\end{align}
\end{proposition}
In Propositions \ref{zx}  and \ref{fd} above, the norms appearing on the right hand sides of those equations are continuous, and of course, independent of the varieties which appear on the left.

The term $R_{G^{\mathbb{C}}|X}$ is explained as follows.
Let
\[ \begin{CD}
\mathfrak{X} @>\pi_{2}>> \cpn \\
@V\pi_{1}VV \\
\mathcal{H}\\
\end{CD}
\]
be a $\gc$ equivariant fibration between smooth varieties. Let $\mathcal{H}_{\infty}$ denote the locus of  smooth fibers in $\mathcal{H}$, i.e. $X_{z}:= \pi_{1}^{-1}(z)$ is a 
smooth subvariety of $\cpn$ for $z \in \mathcal{H}_{\infty}$.
Then $\pi_{2}^{*}\omega_{FS}$ induces a metric on the relative canonical bundle
\begin{align*}
\mathcal{K}_{\mathfrak{X}}:= K_{\mathfrak{X}}\otimes \pi_{1}^{*} K_{\mathcal{H}}^{-1}.
\end{align*} 
over the smooth locus. We let $R_{\mathfrak{X}|\mathcal{H}}$ denote its curvature.
Let $p_{z}$ denote the projection

\[ \begin{CD}
\gc X_{z}@>p_{z}>> \mathfrak{X} @>\pi_{2}>>\cpn\\
@Vp_{1}VV @V\pi_{1}VV \\
\gc@>T_{z}>> \mathcal{H}\\
\end{CD}
\]
Then 
\begin{align*}
R_{\gc |X_{z}}:= p_{z}^{*}R_{\mathfrak{X}|\mathcal{H}}
\end{align*} 
Sometimes we will drop the manifold $X$ in the expression $R_{G^{\mathbb{C}}|X}$ and just write $R_{\gc}$. It will be clear what the variety is from the context.
Observe that $\pi_{2}^{*}\omega_{FS}$ induces a K\"ahler  metric on $\pi_{1}^{-1}(z)$ $(z \in \mathcal{H}_{\infty})$ and hence a metric on the \emph{relative} canonical bundle $K_{X_{z}}$ which we denote by $R(\pi_{2}^{*}(\omega_{FS}))$. 
Now let $g_{\mathfrak{X}}$ and $g_{\mathcal{H}}$ denote two kahler metrics on 
$\mathfrak{X}$ and $\mathcal{H}$ respectively.
In this way we obtain \emph{another} metric on the relative canonical bundle
over the smooth locus. We let $R_{\mathfrak{X}|\mathcal{H}}$ denote its curvature.
\begin{align*}
R_{\mathfrak{X}|\mathcal{H}}:= R(g_{\mathfrak{X}})-\pi_{1}^{*}R(g_{\mathcal{H}})
\end{align*} 

Then there is a smooth function $\Psi$ defined away from $\mathfrak{X}_{sing.}$
\begin{align*}
\Psi:\mathfrak{X}\backslash \pi_{1}^{-1}(\Delta)\rightarrow \mathbb{R}
\end{align*}
such that
\begin{align*}
 R(\pi_{2}^{*}(\omega_{FS}))=R_{\mathfrak{X}}-\pi_{1}^{*}(R_{\mathcal{H}})+\sqrt{-1}\dl\dlb \Psi
\end{align*}
Moreover, if we define $\Psi_{\mathcal{H}}(z):= \int_{\{y\in\pi_{1}^{-1}(z)\}}\Psi(y) \pi_{2}^{*}(\omega_{FS})^{n}$. Then $\Psi_{\mathcal{H}}(z)$ is bounded from above, is smooth outside $\Delta$, continuous outside $\Delta_{m}$, and goes to $-\infty$ as $z\rightarrow \Delta_{{m}}$. $\Delta_{{m}}$ denotes the locus of $z\in \mathcal{H}$ where $X_{z}$ has a component of multiplicity greater than one.
We have the following theorem ([Tian94]).
\begin{theorem}\label{sing}
There is a continuous metric $||\ ||_{CM}^{2}$ on the CM polarisation such that
if $z\in \mathcal{H}_{\infty}$, then 
\begin{align*}
d(n+1)\nu_{\omega,z}(\sigma)= {\log}\left(e^{(n+1)\Psi_{\mathcal{H}}(\sigma z)}\frac{||\ ||_{CM}^{2}(\sigma z)}{||\ ||_{CM}^{2}( z)}\right)
\end{align*}
\end{theorem}
This says that the K-energy map is the logarithm of a singular metric on the CM polarisation.

\section{Proof of Theorem \ref{ext}}
In this section we identify the CM polarisation on $\mathcal{H}$, i.e. we will show that
\begin{align*}
L_{CM}^{-1}= \iota^{*}_{HC}\underline{O}_{{\mbox{Div}}_{(d,\dots,d)}}(-\nu_{1})\otimes {L}_{\Delta}^{\nu_{2}}
\end{align*}
The proof is a straightforward application of the Grothendeicks' very powerful
Reimann Roch Hirzebruch Theorem.
We will constantly appeal to the following diagrams.
\[\begin{CD}
\mathfrak{X}\times_{\cpn}E@>\rho_{2}>>E@>p_{2}>>\gr \\
@V\rho_{1}VV    @Vp_{1}VV \\
\mathfrak{X}@>\pi_{2}>> \cpn \\
@V\pi_{1}VV\\
\mathcal{H}
\end{CD}
\]
In the diagram below, the map $p$ denotes the resolution of singularities, which is given by the obvious projection, $\mathcal{D}_{\mathcal{H}}$ is the singular divisor of the fiber product $\mathcal{H}\times_{B_{d}}\Sigma$, and 
$\Delta$ is the proper transform of this divisor inside $\mathfrak{X}\times_{\cpn}E$. Define $\tau:= \pi_{1}\rho_{1}$. $E$ is the flag manifold constructed in the following way.

Let Q be the canonical quotient bundle on the grassmannian $\gr$  (whose dimension is denoted $l+1$). Let 
$L_{Q}$ denote its determinant, the positive generator of Pic($\gr$).
Often in the text we will use the same notation $L_{Q}$ to mean its pull back under the various maps that appear in the diagrams below.
On the product $\cpn \times \gr$ we have the projections $p_{1}$ and $p_{2}$ onto the first and second factors respectively. Let $U$ be the $\underline{O}(-1)$ sheaf on \cpn . By definition $$ p_{1}^{*}(U)|_{([v],L)} = \Ce v $$ and $$p_{2}^{*}(Q)|_{([v],L)} = \Ce ^{N+1}/ \overline{L}$$ here, $\overline{L}$ is the subspace of $ \Ce^{N+1}$ lying over $L \subset \cpn$ . This gives a natural map:$$\Phi: p_{1}^{*}(U)|_{([v],L)} \ra p_{2}^{*}(Q)|_{([v],L)}$$ defined by taking the quotient by $L$. Therefore $\Phi$ induces a holomorphic section $\overline{\Phi}$ of $p_{1}^{*}(U)^{\vee}\otimes p_{2}^{*}(Q)$. 
Now set 
\begin{align}
E:= \{\overline{\Phi}= 0\} = \{(x,L) \in \cpn \times \gr\ :x \in L\}
\end{align}

\[\begin{CD}
\Delta=p^{-1}(\mathcal{D}_{\mathcal{H}})\subset\mathfrak{X}\times_{\cpn}E\\
@VpVV\\
\mathcal{D}_{\mathcal{H}}\subset\mathcal{H}\times_{B_{d}}\Sigma@>g_{2}>>\Sigma@>f_{2}>>\gr \\
@Vg_{1}VV @Vf_{1}VV \\
\mathcal{H}@>\iota_{HC}>>{\mbox{Div}}_{(d,\dots,d)}\subset B_{d}
\end{CD}
\]

With these notations established, we can begin the argument. By GRR we have
\begin{align*}
c_{1}(L_{CM}^{-1})+d(n+1)c_{1}(\mathcal{H})= {\pi _{1}}_{*}\left(-(n+1)c_{1}(L)^{n}c_{1}(K_{\mathfrak{X}})-\mu c_{1}(L)^{n+1}\right)
\end{align*}
Hence, for any $\eta\in H^{(1,1)}(\mathcal{H})$ we see that
\begin{align*}
\begin{split}
&\int_{\mathcal{H}}\eta\wedge (c_{1}(L_{CM}^{-1})+d(n+1)c_{1}(\mathcal{H}))=\\
&(n+1)\int_{\mathfrak{X}}\pi_{2}^{*}c_{1}(L)^{n}\pi_{1}^{*}\eta - \mu\int_{\frak{X}}\pi_{2}^{*}c_{1}(L)^{n+1}\pi_{1}^{*}\eta = \\
&(n+1)\int_{\mathfrak{X}\times_{\cpn}E}\rho_{2}^{*}p_{2}^{*}c_{1}(L_{Q})^{l}\left(c_{1}(\mathfrak{X}\times_{\cpn}E)-N\rho_{2}^{*}p_{2}^{*}c_{1}(L_{Q}) + (n+1)\rho_{1}^{*}\pi_{2}^{*}c_{1}(L)\right)\wedge \tau^{*}\eta\\
&-\mu\int_{\mathfrak{X}\times_{\cpn}E}\rho_{2}^{*}p_{2}^{*}c_{1}(L_{Q})^{l+1}\tau^{*}\eta 
\end{split}
\end{align*}
Next, we introduce a virtual bundle $\mathcal{E}$ on the space $\Sigma$
\begin{align*}
\mathcal{E}:= (-1)^{l}2(\underline{\mathbb{C}}-L_{Q})^{l}+ (-1)^{l}(\underline{\mathbb{C}}-L_{Q})^{l+1}
\end{align*}
We need the following, for the proof we refer the reader to the foundational paper of Mumford and Knudsen [KnudMum76].
\begin{proposition}(Compatibility of Det with base change)\\
\begin{align*}
\iota_{HC}^{*}{\bf{Det}}_{B_{d}}(R_{{f_{1}}_{*}}^{\bull}\mathcal{E})= {\bf{Det}}_{\mathcal{H}}(R
_{{g_{1}}_{*}}^{\bull}g_{2}^{*}\mathcal{E})
\end{align*}
\end{proposition}
It is not difficult to see that
\begin{align*}
{\bf{Det}}_{B_{d}}(R_{{f_{1}}_{*}}^{\bull}\mathcal{E})=\underline{O}_{B_{d}}(D(N+1-2d))
\end{align*}
Another application of GRR\footnote{Here we must use the intersection theory developed by Fulton and MacPherson, since $\mathcal{H}\times_{B_{d}}\Sigma$ is \emph{singular} in codimension one} yields 
\begin{align*}
\begin{split}
&(n+1)D(N+1-2d)\int_{\mathcal{H}}\iota_{HC}^{*}c_{1}(\underline{O}_{B_{d}}(1))\wedge \eta +dD(n+1)\int_{\mathcal{H}}c_{1}(\mathcal{H})\wedge \eta\\
&= (n+1)\int_{\mathfrak{X}\times_{\cpn}E}c_{1}(L_{Q})^{l}\left(c_{1}(\mathfrak{X}\times_{\cpn}E)
-PD([\Delta])\right)\tau^{*}\eta\\
\end{split}
\end{align*}
Therefore
\begin{align*}
\begin{split}
&D\int_{\mathcal{H}}\eta\wedge (c_{1}(L_{CM}^{-1})- D(n+1)(N+1-2d)\int_{\mathcal{H}}\iota_{HC}^{*}c_{1}(\underline{O}_{B_{d}}(1))\wedge \eta\\
&=D((n+1)^{2}-(N(n+1)+\mu))\int_{\mathcal{H}}\iota_{HC}^{*}c_{1}(\underline{O}_{B_{d}}(1))\wedge \eta +(n+1)\int_{\mathcal{H}}c_{1}(L_{\Delta})\wedge \eta
\end{split}
\end{align*}
Now we can finally conclude (by Kodaira-Serre duality) that
\begin{align*}
c_{1}(L_{CM}^{-1})= \iota_{HC}^{*}c_{1}(\underline{O}_{B_{d}}(-\nu_{1} )) + c_{1}(L_{\Delta}^{\otimes \nu_{2}}) 
\end{align*}
From which the theorem follows.$\Box$

Let us explain the line bundle $L_{\Delta}$.
Just set $\mathcal{E}_{\Delta}$ 
to be the following virtual bundle on $\mathfrak{X}\times_{\cpn}E$
\begin{align*}
\mathcal{E}_{\Delta}:=    (-1)^{l+1}(\underline{\mathbb{C}}-\underline{O}(\Delta))(\underline{\mathbb{C}}-L_{Q})^{l}
\end{align*}
Then
\begin{align*}
       L_{\Delta}:= {\bf{Det}}_{\mathcal{H}}(R_{{\tau}_{*}}^{\bull}\mathcal{E}_{\Delta})
\end{align*}

Now let us compare the above result with Theorem 3.7. 
In order to do this, let $|| \ ||$ be any hermitian metric on $\underline{O}(\Delta)$. Then we get a continuous metric on $L_{\Delta}$ by integration of
$C_{1}(\underline{O}(\Delta))C_{1}(L_{Q})$ over the fiber. Call this metric $|| \ ||_{\Delta}$.
It follows from what we have just done that
\begin{align}
\log\left(\frac{||\ ||_{CM}^{2}(\sigma z)}{||\ ||_{CM}^{2}( z)}\right)=
\nu_{1}\log \frac{||R_{X}\circ \sigma^{-1}||^{2}}{||R_{X}||^{2}} + \nu_{2}\log\frac{||\ ||^{2}_{\Delta}(\sigma z)  }{||\ ||^{2}_{\Delta}(z)}
\end{align}
Let $S_{\Delta}$ be a section of $\underline{O}(\Delta)$ corresponding to the divisor $\Delta$. We define a function $\theta_{\Delta}(z)$ on the compliment of the multiplicity locus by
\begin{align*}
\theta_{\Delta}(z):= \int_{\tau^{-1}(z)}\log||S_{\Delta}||^{2}C_{1}(L_{Q})^{l}
\end{align*}
Then $\theta_{\Delta}(z)$ is bounded above and goes to $-\infty$ as $z\ra y\in \mathcal{H}_{m}$.
As a consequence of the Poincare Lelong equation and arguments in the next paragraph we have
\begin{align}
\log\frac{||\ ||^{2}_{\Delta}(\sigma z)}{||\ ||^{2}_{\Delta}(z)}= -D
\log\frac{||(f_{\mathcal{D}_{z}}\circ \sigma^{-1}) ||^{2}}{||f_{\mathcal{D}_{z}}||^{2}}
+\theta_{\Delta}(\sigma z)
\end{align}
It follows from this and Theorem 3.7 that

\begin{align*}
d\nu_{\omega,z}(\sigma)= \Psi_{\mathcal{H}}(\sigma z)+\frac{1}{D}\theta_{\Delta}(\sigma z) +\frac{\nu_{1}}{n+1}\log \frac{||R_{X}\circ \sigma^{-1}||^{2}}{||R_{X}||^{2}} -\log\frac{||(f_{\mathcal{D}_{z}}\circ \sigma^{-1}) ||^{2}}{||f_{\mathcal{D}_{z}}||^{2}}
\end{align*}

 
Now we give the proof of Theorem \ref{asymp}, we begin with a proof of proposition \ref{fd}.
Lets fix some notation. Given $\mathcal{V}$ we define
\begin{align*}
\Gamma_{\mathcal{V}}:= \{(L,E\subset F)\in {\mathcal{V}}\times \mathbf{Fl}_{2} : L\in P_{EF}\}
\end{align*}

\begin{align*}
Z_{\mathcal{V}}\ \mbox{ is the hypersurface in}\  \mathbf{Fl}_{2}\ \mbox{associated to }\ \mathcal{V}
\end{align*}

\begin{align*}
\begin{split}
\Sigma:=\{([f],\ E \subset F)\in  B_{d_{1},d_{2}}\times \mathbf{Fl}_{2}:f( E \subset F)=0\}
\end{split}
\end{align*}

\begin{align*}
T_{f_{\mathcal{V}}}(\sigma):= [f_{\mathcal{V}}\circ \sigma^{-1}]
\end{align*}

\begin{align*}
\gc Z_{\mathcal{V}}:= \{(\sigma,\ E\subset F):f_{\mathcal{V}}(\sigma^{-1}(E\subset F))=0\} 
\end{align*}
$\mathbf{Fl}_{2}$ is the two step flag $\{E\subset F\}$ of subspaces in $\cpn$
with $dim(E)= k-1$ and $dim(F)=k+1$, $\mathbf{Fl}_{3}$ is the  three step flag,
whose middle subspace has dimension $k$, the outer two lie in $\mathbf{Fl}_{2}$.\\
The reader may find the diagrams below to be of some use.

\[ \begin{CD}
\Gamma_{\mathcal{V}}@>>>\mathbf{Fl}_{3}@>p_{2}>> \mathbf{Fl}_{2}\qquad (*)\\
@V\pi_{1}VV  @V p_{1} VV\\
\mathcal{V}@>>> Gr(k,\cpn)\\
\end{CD}
\]
\\
\\

\[ \begin{CD}
G^{\mathbb{C}}Z_{\mathcal{V}}@>T_{f_{\mathcal{V}}}\times id>> \Sigma@>p_{2}>> \mathbf{Fl}_{2} \qquad (**)\\
@V\pi_{1}VV  @V p_{1} VV\\
\gc@>T_{f_{\mathcal{V}}} >> B_{d_{1},d_{2}}\\
\end{CD}
\]

The standard hermitian form on $\mathbf{C}^{N+1}$ induces a metric $h_{fl_{2}}$ on the line bundle
\begin{align*}
\underline{O}(1,1):= \mathbb{L}_{1}\otimes\mathbb{L}_{2}
\end{align*}
We let $C_{1}(\underline{O}(1,1))$ denote the representative of the $1^{st}$ 
chern class of $\underline{O}(1,1)$  with respect to the metric $h_{fl_{2}}$.
Let $\ell_2$ be the dimension of $\mathbf{Fl}_{2}$.\\
\emph{Proof of proposition }(\ref{fd})\\
We begin with the following
\begin{lemma}\label{degfl2}
\begin{align}
{p_{1}}_{*}\ p_{2}^{*}\ C_{1}(\underline{O}(1,1))^{\ell_2}= \frac{1}{D}deg(\mathbf{Fl}_{2})\omega_{Pl}^{(N-k)(k+1)-1}
\end{align}
Where $deg(\mathbf{Fl}_{2}):= \int_{\mathbf{Fl}_{2}}C_{1}(\underline{O}(1,1))^{\ell_2}$
\end{lemma}
To see this, first observe that ${p_{1}}_{*}\ p_{2}^{*}\ C_{1}(\underline{O}(1,1))^{\ell_2}$ is invariant under the action of the unitary group, since $h^{(1,1)}(Gr(k,\cpn))= 1$ we must have
\begin{align*}
{p_{1}}_{*}\ p_{2}^{*}\ C_{1}(\underline{O}(1,1))^{\ell_2}= C\omega_{Pl}^{(N-k)(k+1)-1}
\end{align*}
For some constant $C$. To find $C$ we just integrate against $\omega_{Pl}$
\begin{align*}
CD= \int_{\mathbf{Fl}_{3}}p_{2}^*C_{1}(\underline{O}(1,1))^{\ell_2}\wedge p_{1}^{*}\omega_{Pl}
\end{align*}
The component of $PD(\mathbf{Fl}_{3})$ in $H^{(N-k)(k+1)-1}(Gr(k,\cpn))$ is given by the poincare dual to a generic slice
\begin{align*}
PD(\mathbf{Fl}_{3}\cap (Gr(k,\cpn)\times \{E\subset F\})) = p_{1}^{*}PD(P_{EF})
\end{align*}
Here, $P_{EF}$ is the corresponding pencil in the grassmannian.
Since
\begin{align*}
\int_{P_{EF}}\omega_{Pl}=1
\end{align*}
we see that
\begin{align*}
PD(P_{EF})= \frac{1}{D}\omega_{Pl}^{(N-k)(k+1)-1}\qquad \Box
\end{align*}
Next, we define a $(1,1)$ current $U_{d_{1},d_{2}}$ on $B_{d_{1},d_{2}}$ as follows
\begin{align*}
\int_{B_{d_{1},d_{2}}}U_{d_{1},d_{2}}\wedge \varphi := \int_{\Sigma}\pi_{1}^{*}\varphi \wedge \pi_{2}^{*}C_{1}(\underline{O}(1,1))^{\ell_2}
\end{align*}
Then we have the following
\begin{proposition}
There is a continuous metric $|| \ ||$ on $\underline{O}_{B_{d_{1},d_{2}}}(-1)$ such that
\begin{align*}
\frac{\sqrt{-1}}{2\pi}{\mbox{deg}(\mathbf{Fl}_{2})}\ \dl\dlb\log|| \ ||^{2} =  U_{d_{1},d_{2}}
\end{align*}
in the sense of distributions.
\end{proposition}

This holds for any family coming from a very ample complete linear system, 
therefore below we will consider the data of a compact complex manifold $M^{n+1}$ and a very ample (hermitian) line bundle $\mathcal{L}^{d}$ of degree $d$.
Its curvature form is denoted by $\omega_{P}$.\\
{\bf{Proof}}\\
Let $\mathfrak{X}$ be the universal family of hypersurfaces of degree $d$ inside of $M$:
\begin{align*}
\mathfrak{X}:= \{([S],p)\in B_{d}\times M : S(p)=0\}
\end{align*}
 
We have the factor projections:
\[ \begin{CD}
{\mathfrak{X}}@>p_{2}>> M\\
@Vp_{1}VV \\
 {B_{d}}\\
\end{CD}
\]
\\
\\
$B_{d}$ is the complete linear system of degree $d$ on $M$
\begin{align*}
B_{d}= \mathbb{P}(H^{0}(\mathcal{L}^{d},M))
\end{align*}
 
Recall that we are studying the $(1,1)$ current $\omega_{B}$ on $B_{d}$:
\begin{align*}
\int_{B_{d}}\varphi\wedge \omega_{B}:= \int_{\mathfrak{X}}p_{1}^{*}\varphi\wedge p_{2}^{*}\omega_{P}^{(n+1)}
\end{align*}
$\varphi$ is a smooth test form of type $(b-1,b-1)$ on $B_{d}$, $b:=dim(B_{d})$. In fact we will prove that $\omega_{B}$ is the curvature form of a \emph{H\"older} continuous metric, although the continuity is all that is really relevant here. First we need an estimate of $\omega_{B}$ near those points $[S]$ where the fiber is singular, notice that $\omega_{B}$ is smooth away from the discriminant locus $\mathfrak{D}$ , which is the hypersurface in $B_{d}$ consisting of singular hypersurfaces in $M$.
\begin{lemma}
Let $\omega_{FS}$ be the K\"ahler form of the Fubini-Study metric on $B_{d}$.
Then there is a constant $C$ such that for any $r<1$,
$$
\int_{B_{r}([S])}\omega_{FS}^{b-1}\wedge\omega_{B}\leq Cr^{2b-2+\frac{2}{d}}
$$
where $B_{r}([S])$ denotes the geodesic ball in $B_{d}$ with radius $r$  
centered at $[S]$.
\end{lemma}
{\bf{Proof}}\\
Choose local coordinates $z_{1},\dots,z_{b}$ at $[S]$ such that $[S]=(0,\dots,0)$, and for any $0<c_{1},\dots,c_{b-1}$, the intersection set $\{z_{i}=c_{i}|i\neq j\}\cap \mathfrak{D}$ is finite, for all $j$. It is enough to show that:
\begin{align*}
\left(\frac{2\pi}{\sqrt{-1}}\right)^{b-1}\int_{B_{r}([S])}\omega_{B}\wedge dz_{1}\wedge d\overline{z}_{1}\wedge \dots\wedge \widehat{dz_{j}}\wedge\widehat{d\overline{z}_{j}}\wedge \dots \wedge dz_{b}\wedge d\overline{z}_{b}\leq Cr^{2(b-1)+\frac{2}{d}}
\end{align*}
where $j=1,2,\dots,b$.
For simplicity let us take the last one, so $j=b$. Define a holomorphic map:
\begin{align*}
f:p_{1}^{-1}(B_{r}([S]))\rightarrow B_{r}(0)\cap \mathbb{C}^{b-1}\\
f((z_{1},\dots,z_{b}),p ):= (z_{1},\dots,z_{b-1})
\end{align*}
Where
\begin{align*}
 B_{r}(0)\cap \mathbb{C}^{b-1}=\{|z_{1}|^{2}+\dots + |z_{b-1}|^{2}<r^{2}\}.
\end{align*}
By definition of $\omega_{B}$ we have,\\
\begin{align*}
\begin{split}
&\left(\frac{2\pi}{\sqrt{-1}}\right)^{b-1}\int_{B_{r}([S])}\omega_{B}\wedge dz_{1}\wedge d\overline{z}_{1}\wedge \dots\wedge {dz_{j}}\wedge{d\overline{z}_{j}}\wedge \dots \wedge dz_{b-1}\wedge d\overline{z}_{b-1}\\
=& \left(\frac{2\pi}{\sqrt{-1}}\right)^{b-1}\int_{p_{1}^{-1}(B_{r}([S]))}f^{*}( dz_{1}\wedge d\overline{z}_{1}\wedge \dots\wedge  {dz_{j}}\wedge {d\overline{z}_{j}}\wedge \dots \wedge dz_{b-1}\wedge d\overline{z}_{b-1})\wedge p_{2}^{*}\omega_{P}^{(n+1)}\\
 &\leq \int_{p_{1}^{-1}(B_{r}([S]))}|\dl f|^{2}\omega^{b+n}
\end{split}
\end{align*}
where $\omega$ is the product of $\omega_{FS}$ and $\omega_{P}$ on $B_{d}\times M$. In this proof $C$ denotes a uniform constant.
By the Co-Area formula:
\begin{align*}
\begin{split}
&\int_{p_{1}^{-1}(B_{r}([S]))}|\dl f|^{2}\omega^{b+n}\\
& = \int_{\{|z_{1}|^{2}+\dots + |z_{b-1}|^{2}<r^{2}\}}Vol(f^{-1}(z_{1},\dots,z_{b-1}))dz_{1}\wedge d\overline{z}_{1}\wedge \dots \wedge dz_{b-1}\wedge d\overline{z}_{b-1}
\end{split}
\end{align*}
Since $|\dl f|\leq \sqrt{b}$ on $p_{1}^{-1}(B_{r}([S]))$, it reduces to show that:
\begin{align*}
Vol(f^{-1}((z_{1},\dots,z_{b-1}))\leq Cr^{\frac{2}{d}}.
\end{align*}

Fix any $(b-1)$ tuple $(c_{1},\dots, c_{b-1})$, define $\mathcal{C}$ to be the holomorphic curve:

\begin{align*}
\mathcal{C}:= \{z_{i}=c_{i}, |z_{b}|<1|1\leq i \leq b-1\}
\end{align*}
we want to estimate $Vol(p_{1}^{-1}(\mathcal{C}\cap B_{r}([S])))$ with $[S]\in \mathfrak{D}$. Recall that $p_{1}^{-1}([S])$ is a hypersurface in $M$. Cover $p_{1}^{-1}([S])$ with finitely many coordinate charts $\{(U_{\alpha};z_{\alpha 1},\dots,z_{\alpha n+1}\})\}$ of $\mathcal{C}\times M$ satisfying:
\begin{flushleft}(1) For some constant $C>0$,
\begin{align*}
\begin{split}
Vol(p_{1}^{-1}(\mathcal{C}\cap B_{r}([S])))\leq & \Sigma_{\alpha}Vol(U_{\alpha}\cap p_{1}^{-1}(\mathcal{C}\cap B_{r}([S])))\\
&\leq CVol(p_{1}^{-1}(\mathcal{C}\cap B_{r}([S])));
\end{split}
\end{align*}
\\
(2) For some positive number $\delta$, $U_{\alpha}:= \{|z_{\alpha 1}|<\delta,\dots 
|z_{\alpha n+1}|<\delta\}$, and for any $a_{1},\dots, a_{n-1}$ with the $|a_{i}|<\delta$, the intersection of the complex line segment $\{z_{\alpha i}=a_{i}|1\leq i \leq n-1\}$ and $p_{1}^{-1}([S])$ consists of less than $d$ points counted with multiplicity.
\end{flushleft}
Consider the restriction
\begin{align*}
p_{1}:\mathcal{S}\rightarrow \mathcal{C}\cap B_{r}([S])
\end{align*}
where 
\begin{align*}
\mathcal{S}:= \{z_{\alpha i}= a_{i}|1\leq i\leq n-1\} \cap p_{1}^{-1}(\mathcal{C}\cap B_{r}([S])).
\end{align*}

This is a branch map of degree not greater than $d$. Then
\begin{align*}
\{z_{\alpha i}=a_{i}:1\leq i \leq n-1\}\cap p_{1}^{-1}(\mathcal{C}\cap B_{r}([S]))
\subset B_{cr^{\frac{1}{d}}}(p_{1}^{-1}([S])\cap U_{\alpha})
\end{align*}
where $c$ is independent of $a_{i}$. It follows that:
\begin{align*}
p_{1}^{-1}(\mathcal{C}\cap B_{r}([S]))\subset B_{cr^{\frac{1}{d}}}(p_{1}^{-1}([S]))
\end{align*}
The required volume estimate is now immediate. This completes the proof of lemma 4.2.
$\Box$

Recall that a function on $\mathbb{R}^{k}$ is subharmonic if its Laplacian in nonnegative in the sense of distributions.
\begin{lemma}
Let $\varphi$ be a subharmonic function on the unit ball $B_{1}(0)$ in 
$\mathbb{R}^{k}$ satisfying:
\begin{flushleft}

$\bullet$  $\int_{\dl B_{1}(0)}|\varphi|dS\leq C$\\

$\bullet$ There is a constant $\epsilon \in (0,\frac{1}{2})$ such that for every
$x\in B_{1}(0)$, and $r\leq dist(x,\dl B_{1}(0))$
$$ \int_{B_{1}(0)\cap B_{r}(x)}\triangle \varphi dV \leq Cr^{k-2 +\epsilon}$$
where $C$ is a constant.\\
\end{flushleft}

Then there is a constant $C^{'}$ depending only on $C$, $\epsilon$, and $k$, such that for every $x,y \in B_{\frac{1}{2}}(0)$,
$$|\varphi(x)- \varphi(y)|\leq C^{'}|x-y|^{\epsilon}.$$

In particular, $\varphi$ is H\"{o}lder continuous on $B_{\frac{1}{2}}(0)$.
\end{lemma}
{\bf{Proof}}
In this proof, we will always use $C^{'}$ to denote a constant depending only on $C$, $\epsilon$, and $k$. As usual, its actual value will vary from line to line.
We just consider the case $x=z$ and $y=0$. We let $K(z,\zeta)$ be the Greens function of the unit ball,
\begin{align*}
K(z,\zeta):= \frac{1}{|z-\zeta|^{k-2}}-\frac{|z|^{k-2}}{|z-|z|^{2}\zeta|^{k-2}}
\end{align*}
By the Green formula:
\begin{align*}
\begin{split}
&\varphi(z)- \varphi(0)= \int_{B_{1}(0)}\triangle \varphi(K(0,\zeta)-K(z,\zeta))dV \\ 
&+\int_{\dl B_{1}(0)}\varphi \left(\frac{\dl K}{\dl \nu}(0,\zeta)-
\frac{\dl K}{\dl \nu}(z,\zeta)\right)dS
\end{split}
\end{align*}

Using our first assumption and the smoothness of $\frac{|z|^{k-2}}{|z-|z|^{2}\zeta|^{k-2}}$ on the ball, and that of $\frac{\dl K}{\dl \nu}(0,\zeta)$ and
$\frac{\dl K}{\dl \nu}(z,\zeta)$ on its boundary, we can deduce that:
\begin{align*}
\begin{split}
&|\varphi(z)- \varphi(0)|\leq C^{'}|z| + |\int_{B_{1}(0)}\left(\frac{1}{|z-\zeta|^{k-2}}-\frac{1}{|\zeta|^{k-2}}\right)\Delta \varphi dV|\\
&\leq C^{'}|z| + |\int_{B_{2|z|}(0)}\left(\frac{1}{|z-\zeta|^{k-2}}-\frac{1}{|\zeta|^{k-2}}\right)\Delta \varphi dV|\\
&+ \sum_{j\geq 1; 2^{j+1}|z|\leq 1}|\int_{B_{2^{j+1}|z|}(0)\backslash B_{2^{j}|z|}(0)}\left(\frac{1}{|z-\zeta|^{k-2}}-\frac{1}{|\zeta|^{k-2}}\right)\Delta \varphi dV|
\end{split}
\end{align*}

It follows from the mean value theorem that, for $|\zeta|\geq 2|z|$,
\begin{align*}
|\frac{1}{|z-\zeta|^{k-2}}-\frac{1}{|\zeta|^{k-2}}|\leq (k-2)2^{2k-4}\frac{|z|}{|\zeta|^{k-1}}.
\end{align*}

Therefore,
\begin{align*}
\begin{split}
&\sum_{j\geq 1; 2^{j+1}|z|\leq 1}|\int_{B_{2^{j+1}|z|}(0)\backslash 
B_{2^{j}|z|}(0)}\left(\frac{1}{|z-\zeta|^{k-2}}-\frac{1}{|\zeta|^{k-2}}\right)\Delta \varphi dV|\\
&\leq (k-2)2^{2k-4}\sum_{j\geq 1; 2^{j+1}|z|\leq 1}2^{-j(k-1)}|z|^{-k+2}\int_{
B_{2^{j+1}|z|}(0)}\Delta \varphi dV\\
&\leq (k-2)2^{2k-4}\sum_{j\geq 1; 2^{j+1}|z|\leq 1}2^{-j(k-1)}|z|^{-k+2}2^{(j+1)(k-2+\epsilon)}|z|^{k-2+\epsilon}\\
&\leq (k-2)2^{3k-6+\epsilon}|z|^{\epsilon}\sum_{j\geq 1}2^{-j(1-\epsilon)}\\
&\leq C^{'}|z|^{\epsilon}.
\end{split}
\end{align*}
It follows from the above that,
\begin{align*}
\begin{split}
&|\varphi(z)- \varphi(0)|\leq C^{'}|z|^{\epsilon}+\int_{B_{|z|}(0)}
\frac{1}{|z-\zeta|^{k-2}}\Delta \varphi dV\\
& + \int_{B_{|z|}(0)}
\frac{1}{|\zeta|^{k-2}}\Delta \varphi dV.
\end{split}
\end{align*}
Now we have,

\begin{align*}
\begin{split}
&\int_{B_{|z|}(0)}
\frac{1}{|z-\zeta|^{k-2}}\Delta \varphi dV\\
&\sum_{j\geq 1}\int_{B_{2^{-j+1}|z|}(0)\backslash B_{2^{-j}|z|}(0)}\frac{1}{|z-\zeta|^{k-2}}\Delta \varphi dV\\
&\leq \sum_{j\geq 1}(2|z|)^{-j(k-2)}\int_{B_{2^{-j+1}|z|}(0)}\Delta \varphi dV\\
& \leq 2^{k-2 +\epsilon}|z|^{\epsilon}\sum_{j\geq 1}2^{-\epsilon j}= C^{'}|z|^{\epsilon}
\end{split}
\end{align*}
This completes the proof of the lemma, combining this with the previous lemma proves proposition 4.1.
$\Box$

Now we continue with the proof of Proposition (3.2).
Applying Fubinis' theorem to the second diagram $(**)$ gives
\begin{align*}
\int_{\gc Z_{\mathcal{V}}}\pi_{2}^{*}C_{1}(\underline{O}(1,1))^{\ell_2}\wedge \pi_{1}^{*}\varphi =
\int_{\gc}\mbox{deg}(\mathbf{Fl}_{2})\dl\dlb \log || \ ||^{2}\wedge \varphi
\end{align*}
for every smooth compactly supported $(\mbox{dim}\gc-1,\mbox{dim}\gc-1)$ form $\varphi$
on $\gc$.
On the other hand, using lemma \ref{degfl2}, and the fact that $\gc Z_{\mathcal{V}}$
and $\gc \Gamma_{\mathcal{V}}$ are birational we have
\begin{align*}
\int_{\gc Z_{\mathcal{V}}}\pi_{2}^{*}C_{1}(\underline{O}(1,1))^{\ell_2}\wedge \pi_{1}^{*}\varphi =
\frac{\mbox{deg}(\mathbf{Fl}_{2})}{D}\int_{\gc \mathcal{V}}p_{2}^{*}\omega_{Pl}^{(N-k)(k+1)-1}\wedge \varphi
\end{align*}
This completes the proof of proposition \ref{fd} $\Box$\\
\\

Now we recall the construction of the Quillen metric.\\
Let constants be given by
\begin{align}
\bull  C_{1}(m)= D \left(\frac{d(N-n)(n+1)-(N+1-d)}{(N-n)(n+1)}\right)
\end{align}

\begin{align}
\bull C_{2}(m)= \frac{d(N-n)(n+1)}{d(N-n)(n+1)-(N+1-d)}
\end{align}

\begin{align}
\bull C_{3}(m)= \frac{d-N-1}{d(N-n)(n+1)-(N+1-d)}
\end{align}
Let $B_{ss}$ be the supersingular (reducible polynomials with multiplicity, e.g. $f_{0}^{2}f_{1}$,\newline $\mbox{deg}f_{0}>0$) locus of $B_{d}$.
On $B_{d} \setminus B_{ss}$ we introduce the $C^{0}$ function
\begin{align} 
\Psi_{B}([f]):= \int_{Z_{f}}\Psi \ p_{2}^{*}\omega_{Pl}^{(N-n)(n+1)-1}
\end{align}
Where $\Psi$ is the function on $\Sigma_{\infty}$ given by
\begin{align}\label{D}
\Psi([f],L):= \log \left(\frac{||D_{h}f||^{2}(L)}{|||f|||^{2}}\right)
\end{align}

$\Sigma$ is the universal hypersurface

\begin{align} 
\Sigma:= \{([f],L):f(L)=0\} \subset B \times Gr(N-n,\mathbf{C}^{N+1})
\end{align}
With factor projections $p_{1}$ and $p_{2}$

\[ \begin{CD}
\Sigma@>p_{2}>> Gr(N-n,\mathbf{C}^{N+1}) \\
@V p_{1}VV \\
B_{d} \\
\end{CD}
\]

In (\ref{D}) $D_{h}$ is  the holomorphic hermitian connection on $O(d)$ induced by $||.||_{O(-1)}$ and $||| \ |||$ is any norm on the vector space of sections of $O(d)$ . For example,
when $N = n+1 $ we see that $X$ is a hypersurface defined by some homogeneous polynomial $f$ (of degree $d$)
\begin{align*}
f(z_{0},\dots,z_{n+1})= \sum_{i_{0}+ \dots + i_{n+1}= d}A_{i_{0}\dots i_{n+1}}z_{0}^{i_{0}}\dots z_{n+1}^{ i_{n+1}}
\end{align*}
In this case there is an explicit formula for $\Psi$
\begin{align}\label{psi}
\Psi([f],\ [z_{0},\dots,z_{n+1}])= \log\left(\frac{\sum_{i=0}^{n+1}|\frac{\dl f}{\dl z_{i}}|^{2}}{\sum_{i_{0}\dots i_{n+1}}|A_{i_{0}\dots i_{n+1}}|^{2}||z||^{2d-2}}\right)
\end{align}
Then $\Psi$ is well defined away from the singular locus and has the property that

\begin{align*} 
-\frac{\sqrt{-1}}{2\pi}\dl\dlb \Psi = Ric(\omega_{P}|_{Z_{f}}) - (N+1-d)\omega_{P}+ p_{1}^{*}\omega_{B}.
\end{align*}

Now we can define the Quillen metric by
\begin{align} 
|| \ ||_{Q}:= e^{(C_{1}(m)^{-1}\frac{\Psi_{B}}{2})}|| \ ||_{\underline{O}(-1)}^{C_{2}(m)}
|| \ ||_{u}^{C_{3}(m)}
\end{align} 

$|| \ ||_{u}$ is a H\"{o}lder continuous metric on $\underline{O}(-1)_{B}$ defined by the current equation (see Prop. (4.1))
\begin{align}
{p_{1}}_{*}{p_{2}}^{*}\omega_{Pl}^{(N-n)(n+1)}= D\frac{\sqrt{-1}}{2\pi}\dl\dlb\log || \ ||_{u}^{2}
\end{align}
$D$ is the degree of the Grassmannian in its Pl\"{u}cker embedding.
Given 
\begin{align}
f \in H^{0}(\gr, \underline{O}(d))
\end{align}
we set
\begin{align}\label{d0}
D_{0}(\sigma):= C_{1}(m)\log \frac{||(f\circ \sigma^{-1}) ||_{Q}^{2}}{||\ f\  ||_{Q}^{2}}
\end{align} 

We also need the following proposition.
Assume  that $Z_{f}$ is a normal irreducible hypersurface in $\gr$, then the same argument as in [Tian94] shows that 
\begin{proposition}\label{hypgr}
\begin{align}\label{D0=int}
\int_{G}D_{0}(\sigma)\dl\dlb\varphi =
\int_{GZ_{f}}(R_{\gc}+\frac{\mu(Z_{f})}{(N-n)(n+1)}p_{2}^{*}(\omega_{P}))\wedge
p_{2}^{*}(\omega_{P})^{(N-n)(n+1)-1}\wedge p_{1}^{*}(\varphi)
\end{align} 
For every smooth compactly supported test form $\varphi$.
\end{proposition}

To proceed with the proof of Theorem \ref{weight} we would like to use proposition \ref{hypgr}. However, $Z_{X}$ has singularities in codimension one so that this proposition does not, at first, apply. 
The idea is simply to let a sequence of smooth $Z_{f_{j}}$ degenerate to the chow point of $X_{z}$ for $z \in \mathcal{H}_{\infty}$, that $X$ moves in a family is irrelevant at the moment, so we drop the subscript $z$. What we need is the following version of the Poincare Lelong formula\\
\emph{ Let $M$ be a compact complex manifold together with a metrised very ample line bundle $L$, with curvature K\"ahler form $\omega $. Assume that  
\begin{align}
f_{j} \in \mathbf{P}H^{0}(M,L)
\end{align} 
is a sequence of smooth hypersurfaces converging to $f_{\infty}$, where $Z_{\infty}$ has a  codimension one
singular divisor $\mathcal{D}$.
Then for any test $(n-1,n-1)$ form $\varphi$ we have
\begin{align}
\lim_{j}\int_{Z_{j}}Ric(\omega |_{Z_{j}})\wedge \varphi = \int_{Z_{\infty}^{reg.}}
Ric(\omega |_{Z_{\infty}^{reg.}})\wedge \varphi - \int_{\mathcal{D}}\varphi
\end{align}  
}
This follows from the explicit formula for the Ricci potential (\ref{psi}).
 
In our situation we just degenerate the right hand side of (\ref{D0=int}). There are two terms to study
\begin{align}\label{term1}
\bull \int_{G^{\mathbb{C}}Z_{j}}R_{G^{\mathbb{C}}|Z_{j}}\wedge
p_{2}^{*}(\omega_{P})^{(N-n)(n+1)-1}\wedge \pi^{*}\varphi
\end{align}

and
\begin{align}
\bull \frac{\mu(Z_{f_{j}})}{(N-n)(n+1)}\int_{G^{\mathbb{C}}Z_{X}}
p_{2}^{*}(\omega_{P})^{(N-n)(n+1)}\wedge \pi^{*} \varphi
\end{align}

We consider  $\mu(Z_{f_{j}})$ first, this will lead to an explicit formula for the volume of
$\mathcal{D}$. We do not need this formula for our present purpose. However, we do need the calculation which leads to it. Recall that  $\mu(Z_{f_{j}})$ is the average of the scalar curvature
\begin{align}
  \frac{\mu(Z_{f_{j}})}{(N-n)(n+1)-1}:= \frac{1}{dD}\int_{Z_{f_{j}}}Ric(\omega_{P}|_{Z_{f_{j}}})\wedge \omega_{P}^{\zto}\equiv 
N+1-d 
\end{align}
Where we used the adjunction formula in the rightmost equality together with the fact that 
\begin{align*}
Ric(\omega_{P}|_{Z_{f_{j}}})= C_{1}(K^{-1}_{Z_{f_{j}}}).
\end{align*}
According to Poincare Lelong, we have\\
\begin{align}
\begin{split}\label{plelong}
N+1-d =& \frac{1}{dD}\int_{Z^{reg.}_{X}}Ric(\omega_{P})\wedge \omega_{P}^{\zto}\\&-\frac{1}{dD}\int_{\mathfrak{D}}\omega_{P}^{\zto}
\end{split}
\end{align} 
The first integral can be lifted to the resolution $\Gamma$ of $Z_{X}$
\begin{align}
\Gamma:=\{(x,L)\in X\times \gr:x\in L\}
\end{align}  

\begin{align}\label{zreg}
\int_{Z^{reg.}_{X}}Ric(\omega_{P})\wedge \omega_{P}^{\zto}
=\int_{\Gamma}Ric(p_{2}^{*}\omega_{P})\wedge p_{2}^{*}\omega_{P}^{\zto} 
\end{align}  
The calculation of the right hand side of \ref{zreg} is straightforward.
Consider the canonical section $\phi$ of the bundle $p_{1}^{*}\mathcal{U}^{\vee}\otimes p_{2}^{*}Q$
\begin{align}
\phi \in H^{0}(\cpn \times \gr,\ \underline{O}(p_{1}^{*}\mathcal{U}^{\vee}\otimes p_{2}^{*}Q))
\end{align}
$\mathcal{U}$ is the universal bundle on $\cpn$, and $Q$ is the quotient bundle on the grassmannian.
We let $E$ denote the base locus of $\phi$ 

\begin{align*}
E:= \{\phi = 0\}= \{(y,L)\in \cpn \times \gr : y\in L\}\quad (\mbox{a flag variety})
\end{align*}
Then we have
\begin{align}
\Gamma= E \cap (X\times \gr)
\end{align} 
 and hence we observe that $\Gamma$ is just the pullback of $E$ to $X$.
From this we get an exact sequence
\[ \begin{CD}
0@>>>T^{(1,0)}_{\Gamma}@>>>T^{(1,0)}_{ X\times \gr}|_{\Gamma}@> \phi_{*}>>p_{1}^{*}\mathcal{U}^{\vee}|_{X}\otimes p_{2}^{*}Q|_{\Gamma}@>>>0
\end{CD}
\]

Using the Bott-Chern forms (see [BGSI88] section f) we see that 
\begin{align}
c_{1}(K_{\Gamma},p_{2}^{*}\omega_{Pl})= \pi_{1}^{*}c_{1}(K_{X},\omega_{FS})-Np_{2}^{*}\omega_{Pl} + (n+1)\pi_{1}^{*}\omega_{FS} + \dl\dlb \theta \quad \theta \in C^{\infty}(E)
\end{align} 

We will also need the following push forward formulas (Lemma \ref{push}). The relevant diagram is
\[ \begin{CD}
\Gamma@>\iota>>E@>p_{2}>> Gr(N-n,\mathbf{C}^{N+1}) \\
@V\pi_{1}VV  @V p_{1} VV\\
X@>>>\cpn\\
\end{CD}
\]

Let $\beta$ be the (positive) constant
\begin{align}
\beta := \int_{\gr} c_{2}(Q)\omega_{Pl}^{\zto}
\end{align}

\begin{lemma}\label{push}
\begin{align}
{p_{1}}_{*}p_{2}^{*}\omega_{Pl}^{\dmz}= D\omega_{FS}^{n}
\end{align}

\begin{align}
{p_{1}}_{*}p_{2}^{*}\omega_{Pl}^{(N-n)(n+1)}= D\omega_{FS}^{n+1} 
\end{align}

\begin{align}
{p_{1}}_{*}p_{2}^{*}\omega_{Pl}^{(N-n)(n+1)-2}= \beta \omega_{FS}^{n-1}
\end{align}
\end{lemma}
{\bf{Proof}}\newline
Note that
\begin{align*}
PD(E)= C_{n+1}(p_{1}^{*}(U)^{\vee}\otimes p_{2}^{*}(Q))=\sum_{i=0}^{n+1}C_{1}(p_{1}^{*}(U)^{\vee})^{n+1-i}C_{i}(  
p_{2}^{*}(Q))
\end{align*}
now proceed as in lemma \ref{degfl2}
${}\Box{}$\\
This gives us the following

\begin{align}
\begin{split}
\int_{\Gamma}Ric(p_{2}^{*}\omega_{P})\wedge p_{2}^{*}\omega_{P}^{\zto} 
=\  NDd + (\frac{\mu(X)}{n} -(n+1))d\beta
\end{split}
\end{align}

Plugging this into (\ref{plelong}) gives the volume.
\begin{proposition}(Volume of the singular divisor)\\
\begin{align}
\begin{split}
 Vol(\mathcal{D})=  d(d-1)D + \frac{\mu(X)-n(n+1)}{n}d\beta
\end{split}
\end{align}
\end{proposition}

Recall that $G^{\mathbb{C}}\Gamma$ is just the fiber product
\[ \begin{CD}
G^{\mathbb{C}}\Gamma@>\pi_{2}>>E \\
@V\pi_{1}VV  @V p_{1} VV\\
G^{\mathbb{C}}X@>>>\cpn\\
\end{CD}
\]

Then the integral (\ref{term1}) can be handled in exactly the same way. Pulling back to $G^{\mathbb{C}}\Gamma$ gives
\begin{align}
R_{\gc|\Gamma}-\pi_{1}^{*}R_{\gc|X}= -Np_{2}^{*}\omega_{Pl} + (n+1)\pi_{1}^{*}\omega_{FS} + \dl\dlb\theta \circ \pi_{2}
\end{align}
Recall that $R_{\gc}$ denotes the curvature of the relative \emph{canonical} bundle.
Now apply Propositions (\ref{zx}) through (\ref{ken}), and the continuity of the Quillen norm away from the supersingular locus of $B_{d}$
\begin{align*}
&\lim_{j}\int_{G^{\mathbb{C}}Z_{j}}R_{G^{\mathbb{C}}|Z_{j}}\wedge
p_{2}^{*}(\omega_{P})^{(N-n)(n+1)-1}\wedge \pi^{*}\varphi =\\
&\int_{G^{\mathbb{C}}Z_{X}^{reg.}}R_{G^{\mathbb{C}}Z}\wedge p_{2}\omega_{Pl}^{\dmz}\wedge \pi^{*}\varphi + \int_{G^{\mathbb{C}}\mathcal{D}}p_{2}^{*}(\omega_{Pl}^{\dmz})\wedge \pi^{*}\varphi 
=\\ &D\{d\int_{G^{\mathbb{C}}}\dl\dlb \nu_{\omega}\wedge \varphi
 -\frac{\mu(X)}{n+1}\int_{G^{\mathbb{C}}}\dl\dlb\log{||R_{X}\circ \sigma^{-1} ||_{u}}^{2} \wedge \varphi \}\\
&-\{N -n-1\}D\int_{G^{\mathbb{C}}}\dl\dlb\log{||R_{X}\circ \sigma^{-1} ||_{u}}^{2} \wedge \varphi \\
&+D\int_{G^{\mathbb{C}}}\dl\dlb \log||(f_{\mathcal{D}}\circ \sigma^{-1}) ||^{2}\wedge \varphi
\end{align*}
Since $\theta$ is just a bounded term, we absorbed it into (any) one of the norms, which one is not important.
Therefore we get
\begin{proposition}\label{kdouble}
The function\footnote{Recall that the norm $D_{0}(\sigma)$ (\ref{d0}) is \emph{singular}, with the singularities explicitly given by the Ricci potential. The other norms are all continuous. See also [PhSt02] for a similar formula.}
\begin{align}\label{doub}
-d\nu_{\omega}(\varphi_{\sigma})+ (2d +\frac{\mu(X)}{n+1}-(n+2))\frac{\log||R_{X}\circ \sigma^{-1}||^{2}}{||R_{X}||^{2}} - \frac{\log||(f_{\mathcal{D}}\circ \sigma^{-1}) ||^{2}}{||f_{\mathcal{D}}||^{2}}+ \frac{1}{D}\Psi_{B}(\sigma)
\end{align}


{is pluriharmonic}.
\end{proposition}
Since $\pi_{1}(\gc)=1$ this function is of the form $\log|\xi(\sigma)|^{2}$ where $\xi$ is an entire function. Since $\gc$ has a normal compactification
\begin{align}
\overline{\gc}:= \{((z_{ij});\ w): \mbox{det}(z_{ij})= w^{N+1}\}\subset \mathbb{C}P^{N^{2}+2N}
\end{align}
with a single irreducible divisor at $\infty$, and $\xi$  
is of polynomial growth near this divisor (Proposition (3.1)), we can conclude that the function $\xi$
is actually a constant, of modulos one (as can be seen by plugging in $\sigma = \mbox{identity}$). Therefore the expression (\ref{doub}) is equal to 0. This gives the promised generalisation of Theorem \ref{t94}. 
 This completes the proof Theorem \ref{asymp}. 
${}\Box{}$ \\
We observe that the singular terms are equivalent on \gc
\begin{align}
\frac{1}{D}\Psi_{B}(\sigma z)\approx \Psi_{\mathcal{H}}(\sigma z)+\frac{1}{D}\theta_{\Delta}(\sigma z)
\end{align}
\section{K-Energy Asymptotics and generalised Futaki invariants}
Now we would like discuss the generalised Futaki invariant of Ding and Tian, and the asymptotics of the Mabuchi K-Energy.
Let $X$ be a fano manifold, assumed to be plurianticanonically embedded into $\cpn$.
Let $\lambda_{t}$ be an algebraic 1psg. of $\slnc$, diagonalised on the standard basis of
$\mathbb{C}^{N+1}$
\begin{align}\label{lam} 
\lambda_{t}=
\begin{pmatrix}t^{m_{0}} & 0 & \dots & 0 \\
         0& t^{m_{1}}& 0 & \dots \\
         \dots& \dots & \dots & \dots \\
         0& \dots & \dots & t^{m_{N}}\\ \end{pmatrix}
\end{align}
Where the $m_{i}$ are integers ordered so that $m_{0}>m_{1}>\dots > m_{N}$ and
$\sum_{i=0}^{N}m_{i} = 0$.
We consider the action of $\lambda_{t}$ on $\cpn$ and hence on $X$, let $X_{t}$ denote
$\lambda_{t}X$. The generalised Futaki invariant is concerned with the limit of this variety as
$t\rightarrow 0$. This limit always exists in $\cpn$ as an algebraic cycle and we denote it by $X_{\infty}$, we say that $X$ 
\emph{jumps} to $X_{\infty}$.  


In [DT], it is assumed that $X_{\infty}$ is normal. In this situation 
there is a smooth function
$f$ on $X_{\infty}^{reg.}$ such that
\begin{align}
Ric(\omega_{FS})- \omega_{FS}= \frac{\sqrt{-1}}{2\pi}\dl\dlb f \quad \mbox{on} \ X_{\infty}^{reg.}
\end{align}
Then the generalised Futaki invariant  is given by 
\begin{align}
F_{X_{\infty}}(v):= \int_{X_{\infty}^{reg.}}v(f)\omega^{n}
\end{align}
This is a holomorphic invariant [DT]. We remark that the normality hypothesis insures that the integral on the right hand side of (5.3) is well defined. In the smooth case we just get the classical Futaki invariant. Our work here allows us to enlarge the class of 1psg.'s $\lambda_{t}$ used to give the degeneration, in particular the limit need \emph{not} be normal. Let  $\lambda_{t}$ be such that $X_{\infty}$ is multiplicity free. Let $\phi_{t}$ be defined by
\begin{align*}
\lambda_{t}^{*}\omega_{FS}= \omega_{FS}+\frac{\sqrt{-1}}{2\pi}\dl\dlb \phi_{t}
\end{align*}
Then we can introduce our refinement of the generalised Futaki invariant.
\begin{align}\label{limk}
Re({\mathfrak{F}_{X_{\infty}}}(v)):= \frac{2}{d}\left(\left(2d +\frac{\mu(X)}{n+1}-(n+2)\right)\mu(\lambda,\ R_{X})- \mu(\lambda,\ f_{\mathcal{D}})\right)
\end{align}
Combining what we have done here with arguments in [DT], we have the following
\begin{theorem}
Assume $\lambda_{t}$ is a normal degeneration, then
\begin{align*}
Re({F_{X_{\infty}}}(v))= Re({\mathfrak{F}_{X_{\infty}}}(v))
\end{align*}
\end{theorem}

We can reformulate our results as follows:\\

{\bf{Asymptotics of the Mabuchi Energy}}\\
\emph{For any degeneration $\lambda_{t}$ without multiple fibers we have the following asymptotic behavior as $t\rightarrow 0$}
\begin{align*}
\nu_{\omega_{FS}}(\phi_{t})= -\frac{2}{d}\left(\left(2d +\frac{\mu(X)}{n+1}-(n+2)\right)\mu(\lambda,\ R_{X})- \mu(\lambda,\ f_{\mathcal{D}})\right)log(t)+ O(1)
\end{align*}


When $X$ is a smooth hypersurface (in which case there is no $\mathcal{D}$) 
Lu has examined the right hand side of (\ref{limk}) for an \emph{arbitrary}
$\lambda(t)$. Let $X=\{f=0\}$, then
the set up in [Lu] is as follows.

\begin{align*}
f(z_{0},\dots,z_{n+1})= \sum_{i_{0}+ \dots + i_{n+1}= d}A_{i_{0}\dots i_{n+1}}z_{0}^{i_{0}}\dots z_{n+1}^{ i_{n+1}}
\end{align*}
and let $\lambda_{t}$ be given by $(m_{0},m_{1},\dots,m_{n+1})$ as in \ref{lam}.
Then the slope is given by 
\begin{align}
\mu(\lambda,\ f)= \mbox{Max}\left(m_{0}i_{0}+\dots+ m_{n+1}i_{n+1}:A_{i_{0}\dots i_{n+1}}\neq 0\right)
\end{align}
Now define the piecewise linear function
\begin{align}
\psi(x_{0},x_{1},\dots,x_{n+1})= \mbox{Min}_{i_{0}+ \dots + i_{n+1}= d}\left(-\sum_{k=0}^{n+1}m_{k}i_{k} + \sum_{k=0}^{n+1}i_{k}x_{k}\right)
\end{align}
Let
\begin{align}
\psi_{i}(r)= \psi(0,\dots,r,\dots,0) \ r \ \mbox{in the}\ i^{th}\ \mbox{spot}
\end{align}
Then Lu proves the following
\begin{theorem}([Lu01])
Let $\nu_{\omega_{FS}}(\phi_{t})$ be the K-Energy on a smooth hypersurface in $\mathbb{C}P^{n+1}$. Then
for a generic 1psg. $\lambda_{t}$ we have 
\begin{align*}
-\lim_{t\rightarrow 0}t\frac{\dl}{\dl t}\nu_{\omega_{FS}}(\phi_{t})= \frac{2}{d}\left(\frac{\mu(\lambda,\ f)
(n+2)(d-1)}{n+1}-\sum_{i=0}^{n+1}\int_{0}^{\infty}\frac{\dl}{\dl r}\psi_{i}(r)(\frac{\dl}{\dl r}\psi_{i}(r)-1)dr\right)
\end{align*}
\end{theorem}
This is compatible with what we have done here, since for a multiplicity free degeneration the second term on the right vanishes. The same analysis can be carried out for the term $\Psi_{B}$. This will be taken up in a subsequent paper.

\providecommand{\bysame}{\leavevmode\hbox to3em{\hrulefill}\thinspace}

 \small{\textsc{Department of Mathematics, Columbia University  NY. NY. 10027 }}\\
 {\em{E-mail:}} 
 {\texttt{stpaul@math.columbia.edu}}
 \\
 \small{\textsc{Department of Mathematics, MIT,
 Cambridge, MA 02139}}\\
 {\em{E-mail:}} 
 {\texttt{tian@math.mit.edu}}

\begin{thebibliography}{CGY99}

\bibitem[BGSI88]{Bismut Gillet Soule}
    J.M.Bismut, H.Gillet, C.Soule,
    \emph{Analytic torsion and holomorphic determinant bundles I}
     Comm. Math. Phys.\textbf{115}~ (1988)no.1, 49-78.
    
\bibitem[BGSII88]{Bismut Gillet Soule}
    J.M.Bismut, H.Gillet, C.Soule,
\emph{Analytic torsion and holomorphic determinant bundles II}
     Comm. Math. Phys.\textbf{115}~ (1988)no.1, 79-126.
 
\bibitem[BGSIII88]{Bismut Gillet Soule}
     J.M.Bismut, H.Gillet, C.Soule,
\emph{Analytic torsion and holomorphic determinant bundles III}
     Comm. Math. Phys.\textbf{115}~ (1988)no.1, 301-351.
 
\bibitem[Don01]{Donaldson}
   S. Donaldson,
   \emph{Scalar curvature and projective embeddings I}
   Jour. Diff. Geometry \textbf{59}~(2001)479-522. 
\bibitem[Don02]{Donaldson}
   S. Donaldson,
   \emph{Scalar curvature and Stability of Toric varieties}
    Jour. Differential Geometry \textbf{62}~ (2002) 289-349 
\bibitem[GKZ94]{IM Gelfand, MM Kapranov, AV Zelevinsky}
   IM Gelfand, MM Kapranov, AV Zelevinsky
\emph{Discriminants, Resultants and Multidimensional
Determinants}
Birkhauser, Boston, 1994.
\bibitem[Gies82]{Gieseker}
   D. Gieseker,
   \emph{Lectures on Moduli of Curves}
   Tata Institute of Fundamental Research, Springer, Berlin (1982).

\bibitem[Gies77]{Gieseker}
   D.Gieseker,
   \emph{Global Moduli for Surfaces of General Type}
   Invent. Math.\textbf{43}~(1977), 233--282.
\bibitem[KnudMum76]{KnudMum}
   F.Knudsen-D.Mumford, 
   \emph{The projectivity of the moduli space of stable curves
         I, preliminaries on ``det" and ``Div}
      Math. Scand. 39-1 (1976), 19-55. 
\bibitem[Kol95]{Kollar}
   J.Kollar,
\emph{Rational Curves on Algebraic Varieties}
  December 1995, Springer Verlag 
\bibitem[Lu01]{Lu}
Z.Lu,
\emph{K Energy and K stability on Hypersurfaces}
 \textbf{math.DG/0108009} 
\bibitem[Luo98]{Luo}
   H. Luo,
   \emph{Geometric Criterion for Gieseker-Mumford Stability of Polarized
Manifolds}\
   J. Differential Geometry \textbf{49}~(1998), 577--599.

\bibitem[Mum82]{Mumford}
   D. Mumford,
   \emph{Geometric Invariant Theory}
   2nd enlarged ed., 1982, D. Mumford, J.Fogarty, Ergeb. Math. Grenzeb.
   \textbf{3},Springer, Berlin, Vol. 34, 1982.

\bibitem[Mum77]{Mumford}
   D. Mumford,
   \emph{Stability of Projective Varieties} L'Enseignment Math., $II^{e}$ Serie    Tome XXIII Fascicule \textbf{1-2},1977.

\bibitem[Paul04]{Paul}
    S.Paul,
    \emph{Geometric Analysis of Chow Mumford Stability}\newline
    Advances in Mathematics \textbf{182}~ Issue 2,333-356 (2004)
\bibitem[PhSt02]{D.H.Phong,J.Sturm}
    Phong, Sturm,
    \emph{Stability, energy functionals, and K\"ahler-Einstein metrics}
    \textbf{math.DG/0203254}
\bibitem[Tian97]{Tian}
   G. Tian,
   \emph{K\"{a}hler Einstein metrics with positive scalar curvature}\newline
   Invent. Math. \textbf{130}~(1997), 1--37.

\bibitem[Tian94]{Tian}
   G.Tian,
   \emph{The K-Energy on hypersurfaces and stability} \newline
   Comm. Anal. Geom. \textbf{2}~(1994), 239--265. 
\bibitem[Tian00]{Tian}
   G.Tian,
   \emph{Bott-Chern forms and geometric stability}\newline
   Discrete and Continuous Dynam. Syst. \textbf{6}~(2000),211-220.
\bibitem[Tian02]{Tian}
   G.Tian,
   \emph{Extremal Metrics and geometric stability}\newline
   Special issue for S.S. Chern, Houston J. of Math., \textbf{128} (2002), 432-441.

\bibitem[Vieh89]{Viehweg}
   E.Viehweg,
   \emph{Weak Positivity and Stability of certain Hilbert points}\newline
    Inv. Math., Part I,\textbf{96}~ (1989) 639-667; Part II, 101 (1990), 
    191-223.
\bibitem[Wang02]{X.Wang}
X.Wang,
\emph{Balance point and stability for vector bundles}
Thesis Brandies University.
\bibitem[WenMu99]{Wendland, Mueller}
K.Wendland, W.Mueller,
\emph{Extremal K\"{a}hler Metrics and Ray singer Analytic Torsion}
\textbf{math.DG/9904048}


\bibitem[Zhang96]{Zhang}
   S.Zhang,
   \emph{Heights and reductions of semistable varieties}\newline
    Compositio Math. \textbf{104}~(1996) no. 1, 77-105.
    





\end{thebibliography}
\end{document}